%% file: main.tex
\documentclass[]{interact}

\usepackage{natbib}
\bibpunct[, ]{(}{)}{;}{a}{}{,}

\input{Definitions}

\begin{document}




\title{Digital Annealer for quadratic unconstrained binary optimization: a comparative performance analysis 
}

\author{
\name{Oylum \c{S}eker\textsuperscript{*}, Neda Tanoumand, and Merve Bodur\thanks{*Corresponding author. Email: \href{mailto:oylum.seker@utoronto.ca}{oylum.seker@utoronto.ca}}}
\affil{Department of Mechanical and Industrial Engineering, University of Toronto, Toronto, ON, Canada}
}

\maketitle

\begin{abstract}
Digital Annealer (DA) is a computer architecture designed for tackling combinatorial optimization problems formulated as quadratic unconstrained binary optimization (QUBO) models. 
In this paper, we present the results of an extensive computational study to evaluate the performance of DA in a systematic way in comparison to multiple state-of-the-art solvers for different problem classes. 
We examine pure QUBO models, as well as QUBO reformulations of three  constrained problems,
namely quadratic assignment, quadratic cycle partition, and selective graph coloring, with the last two being new applications for DA.
For the selective graph coloring problem, we also present a size reduction heuristic that significantly increases the number of eligible instances for DA. 
Our experimental results show that despite being in its development stage, DA can provide high-quality solutions quickly and in that regard rivals the state of the art, particularly for large instances. Moreover, as opposed to established solvers, 
within its limit on the number of decision variables, DA's solution times are not affected by the increase in instance size. These findings illustrate that DA has the potential to become a successful technology in tackling combinatorial optimization problems.

\end{abstract}


\input{Introduction}

\input{OperatingPrinciples}
\input{PerformanceAnalysis}

\input{conclusion}  

\section*{Acknowledegments}
The authors would like to thank Fujitsu Laboratories Ltd. and Fujitsu Consulting (Canada) Inc. for providing financial support and access to Digital Annealer at University of Toronto. The authors would also like to thank Hamed Pouya for useful discussions at early stages of this work.


%
%

\bibliographystyle{tfcad}
\bibliography{Abbrev,references}   

\end{document}

%% file: Definitions.tex
\usepackage{amsmath}
\usepackage{amssymb}
\usepackage{url}
\usepackage[colorlinks = true,
            linkcolor = blue,
            urlcolor  = blue,
            citecolor = blue,
            anchorcolor = blue]{hyperref}
\usepackage[utf8]{inputenc}
\usepackage[english]{babel}
\usepackage[table,xcdraw]{xcolor}
\usepackage{booktabs}
\usepackage[labelfont=bf,labelsep=space]{caption}
\usepackage[labelfont=bf]{subcaption}
\usepackage{comment}
\usepackage{tikz}
\usepackage{multirow}
\usepackage{graphicx}
\usepackage{longtable}
\usepackage{siunitx}

\newcolumntype{M}[1]{>{\centering\arraybackslash}m{#1}}

\usepackage{mathtools}

\allowdisplaybreaks

\RequirePackage{fix-cm}

\usepackage{titlesec}
\titleclass{\subsubsubsection}{straight}[\subsection]

\newcounter{subsubsubsection}[subsubsection]
\renewcommand\thesubsubsubsection{\thesubsubsection.\arabic{subsubsubsection}}
\renewcommand\theparagraph{\thesubsubsubsection.\arabic{paragraph}} 

\titleformat{\subsubsubsection}
  {\normalfont\normalsize\itshape}{\thesubsubsubsection}{1em}{}
\titlespacing*{\subsubsubsection}
{0pt}{3.25ex plus 1ex minus .2ex}{1.5ex plus .2ex}

\makeatletter
\renewcommand\paragraph{\@startsection{paragraph}{5}{\z@}%
  {3.25ex \@plus1ex \@minus.2ex}%
  {-1em}%
  {\normalfont\normalsize\bfseries}}
\renewcommand\subparagraph{\@startsection{subparagraph}{6}{\parindent}%
  {3.25ex \@plus1ex \@minus .2ex}%
  {-1em}%
  {\normalfont\normalsize\bfseries}}
\def\toclevel@subsubsubsection{4}
\def\toclevel@paragraph{5}
\def\toclevel@paragraph{6}
\def\l@subsubsubsection{\@dottedtocline{4}{7em}{4em}}
\def\l@paragraph{\@dottedtocline{5}{10em}{5em}}
\def\l@subparagraph{\@dottedtocline{6}{14em}{6em}}
\makeatother

\makeatletter
\newenvironment{taggedsubequations}[1]
 {%
  \addtocounter{equation}{-1}%
  \begin{subequations}%
  \def\@currentlabel{#1}%
  \renewcommand{\theequation}{#1.\arabic{equation}}%
 }
 {\end{subequations}}
\makeatother 

\newcommand{\R}{\mathbb{R}}

\newcommand{\bsubeq}{\begin{subequations}}
\newcommand{\esubeq}{\end{subequations}}
\newcommand{\BI}{\begin{itemize}}
\newcommand{\EI}{\end{itemize}}

\newcommand{\be}{\begin{enumerate}}
\newcommand{\ee}{\end{enumerate}}

\newcommand{\cred}{\color{red!85!black}}

\newcommand{\cblue}{\color{blue!50!black}}

\newcommand{\cfujitsu}{\color{black}}

\newcommand{\nt}[1]{\textbf{\textcolor{blue!75}{#1}}}

\newcommand{\os}[1]{\textbf{\textcolor{green!70!black}{#1}}}

\usepackage{algorithmic}
\usepackage{algorithm}
\usepackage{fullpage}
\usepackage{enumerate}
\usepackage[normalem]{ulem}

\newcommand{\E}{\mathbb{E}} 

\newcommand{\penaltyCoef}{\lambda}

\newcommand{\numColors}{c}

\newcommand{\largeVal}{$>>$}
\newcommand{\partitlestyle}{\bfseries \itshape}

%% file: Introduction.tex
\section{Introduction}
Combinatorial optimization problems (COPs) aim to choose the best solution among a finite set of possible solutions in terms of a given objective function. They constitute an important problem class having diverse applications, such as in scheduling, resource allocation, production planning, transportation, telecommunications and network design \citep{paschos2013applications}. 
As such, there is a vast literature on solution methods for these optimization problems, which can be classified as general-purpose methods and problem-specific methods, or, as exact methods and heuristics. In this paper, we focus on an emerging general heuristic approach to solving COPs. More specifically, we employ the Digital Annealer, a hardware architecture specially designed for solving COPs,  using the so-called Ising formulations,  to generate high-speed and high-quality solutions. We conduct a computational study to assess the performance of this recent approach compared to traditional ones. In addition, we provide two new applications for this promising technology. 
In what follows, we first review the solution methods for COPs and the literature on the Digital Annealer, then provide the details of our study. 

\subsection{Solution methods for COPs}
General-purpose methods are designed to be adapted to every COP. The most commonly used approach is to formulate the problem as an integer program (IP), 
and solve it (exactly)
via a state-of-the-art solver, e.g., GUROBI \citep{optimization2018}, CPLEX \citep{cplex2009v12} and  SCIP \citep{Achterberg2009}. The IP solvers are powerful in the sense that they can usually find optimal solutions to moderate-size COPs in a reasonable amount of time. However, in large-scale problems, more specifically real industrial problems, 
IP solvers are unable to provide even good quality solutions, let alone optimal ones, quickly. This may happen due to the complexity of the problem at hand, time limitation to provide a solution, and size of the instances. This brings the necessity of developing  general-purpose heuristic methods.  Moreover, as also mentioned in \citep{methmod}, in many large-scale cases including most real-life applications, it is sufficient to obtain a good-quality solution. In practice, it is mostly desired to find high-quality feasible solutions in a short amount to be implemented in real time. Other than the need of obtaining solutions for a dynamic decision-making framework, fast generation of feasible solutions can be useful to evaluate potential future scenarios for the system and pick the best one to be implemented. Lastly, the extra effort of finding an optimal solution may not be worthwhile due to potential estimation errors in  many problem parameters. In order to produce good-quality solutions for large-scale COPs, many heuristics \citep{heu1,heu3,heu2} and meta-heuristics \citep{meta1,meta2,meta3}, specifically the ones based on simulated annealing \citep{SA3,SA1,SA2}, are proposed in the literature. 

There are problem-specific methods developed to address particular families of COPs. For instance, the shortest path problem is one of the famous COPs for which efficient (exact) solution methods are proposed, e.g., the Dijkstra's algorithm \citep{dij2,dij1} and label correcting algorithms \citep{labelcor}. Such methods are very powerful in the sense that they can provide high-quality solutions for large-scale problems in a short amount of time. However, 
problem-specific methods exist only for a small portion of COPs.

A recently arising theme in the literature is to formulate a COP as a quadratic unconstrained binary optimization (QUBO) model \citep{qubo} and employ a specialized hardware to solve it. The goal is to propose a general methodology in order to tackle COPs and provide high-quality solutions very quickly, e.g., in a matter of minutes. Although COPs are usually modeled as constrained binary linear programs, solving their QUBO reformulations via a hardware specially designed for QUBOs with the aforementioned aim can indeed outperform solving their natural formulations via an IP solver. Moreover, QUBO problems naturally arise in many applications such as routing \citep{qubo_routing1,qubo_routing2}, scheduling \citep{qubo_sch1,qubo_sch2}, finance \citep{qubo_fin2,qubo_fin1}, and biology \citep{qubo_bio2,qubo_bio1}. These problems are more challenging, thus amenable to these new technologies.

Moore's law was the golden rule that set the pace for modern digital revolution for the last fifty years. By approaching to the end of this era, it became clear that the traditional silicon-based processors cannot effort the current computational burden. Therefore, there is need for new technologies which can handle optimization problems using reasonable time and energy. Based on these observations, new architecture designs have been making a frequent appearance in the digital industry. One of the recent developments is the introduction of quantum annealing-based hardware which is commercialized by D-Wave Systems Inc. Although it is guaranteed that quantum annealing provides high-quality solutions under ideal conditions, this technology in its current development stage can be very expensive and difficult to run \citep{sao2019application}. In order to overcome the disadvantages of quantum annealing-based hardware, Fujitsu designed and developed a new computer architecture implemented in digital circuit, the Digital Annealer (DA), to address complex COPs efficiently  \citep{DAU2}. 


A natural way to use a QUBO approach in solving COPs is to formulate the entire problem under consideration as a QUBO model, and feed it to a solver capable of tackling such models, a promising example to which is the DA.
One can also utilize QUBO modeling approach at different levels of a solution strategy. 
For the quadratic assignment problem for instance, \cite{liu2019modeling} present a local search based iterative method, which is an extension of the so-called {quantum local search} (QLS) algorithm.
While in QLS a QUBO model is constructed for a single subset of selected variables at each step of the local search to find an improving solution, the authors extend it by choosing multiple subsets of variables and solving a separate QUBO model for each at every step.
Other problems addressed in a similar manner include graph partitioning and network community detection problems \citep{shaydulin2019hybrid,shaydulin2019network,ushijimamwesigwa2020multilevel}.
Our QUBO approach in this work is to reformulate the full binary constrained formulations of the problems as QUBO models and solving them with DA and selected established solvers.

\subsection{Applications of DA}

In this section, we briefly review some works from the literature that focus on the use of DA in different application areas.
\cite{physics2019} use DA to solve the so-called spin-glass problem in physics, where the problem is formulated as an Ising model and performance of DA is compared with the simulated annealing and parallel tempering Monte Carlo methods. 
The results show that DA is orders of magnitude faster than the other two methods in solving fully connected spin-glass problems; however, it fails to achieve a speedup for sparse two-dimensional spin-glass problems.

\cite{rahman2019} consider the outlier rejection problem in the context of WiFi-based positioning and represent it as a QUBO model. 
They show that solving the problem with DA significantly improves localization accuracy when compared to state-of-the-art methods. 
Graph problems comprise another application area that has been addressed using DA. 
\cite{javad2019digitally} use DA in order to a solve sensor placement problem by modeling it as a minimum vertex cover problem, and report that DA is able to find optimal solutions to all tested benchmark instances and delivers high-quality solutions for a large suite of randomly generated graphs in comparison to well-known alternatives from the literature. 
Another such study is by \cite{naghsh2019digitally}, in which the authors consider the maximum weighted clique problem on both benchmark and randomly generated instances and report that DA outperforms various existing algorithms as well as the state-of-the-art solver CPLEX.

A real-life application of DA is proposed by \cite{control2019}, where the problem of controlling a large number of automated guided vehicles in a factory in Japan is formulated as a QUBO model. 
The authors solve the resulting QUBO formulation using D-Wave, DA and GUROBI, and compare their results. 
By comparing current versions of D-Wave and DA, they show that the probability of attaining an optimal solution using DA is much higher than that of D-Wave, and that DA
can solve the considered problem instances without any need to divide the original problem into smaller subproblems. 

DA is also used to handle a machine learning problem. 
\cite{IsingReg} propose a novel regularization method based on an Ising model for the training phase of a deep neural network. 
One of the recent works focusing on applications and performance analysis of DA is by \cite{DAU2}, which evaluates the efficiency of DA on the maximum/minimum cut problem and the quadratic assignment problem in comparison to CPLEX. 
The authors report that DA outperforms CPLEX in terms of computation time and solution quality. 
Another recent work is by \cite{cseker2020routing}, where the authors address a routing and wavelength assignment problem with a QUBO approach, and show that solving  it using DA outperforms or is comparable to established techniques coupled with the state-of-the-art solvers.

\subsection{Our contributions}
In this paper, we present the results of an extensive computational study to assess the quality of solutions that DA yields, which, to the best of our knowledge, has not been analyzed in detail in the literature. We first consider problems that are originally formulated as QUBO models. We conduct experiments on a large set of instances from existing QUBO libraries, comparing the performance of DA with three state-of-the-art exact solvers. 
In addition, we analyze the impact of variable ordering on the obtained solution quality for DA.

In order to develop a better understanding of DA's efficiency in different classes of problems, we also consider problems whose natural formulations include constraints, and conduct performance analyses for solving their QUBO reformulations. 
In that regard, we use 
the quadratic assignment problem (QAP) 
since it also has available benchmark instances. 
Moreover, we  
introduce two new applications for DA, the quadratic cycle partition problem (QCPP), and the selective graph coloring problem (Sel-Col), which is a generalization of the classical graph coloring problem. For Sel-Col, we design a heuristic method that significantly reduces problem size, hence greatly increasing the number of instances that can be solved by DA. 

Our contributions are as follows:

\begin{enumerate}
    \item 
    To the best of our knowledge, our study constitutes the first work assessing the performance of DA in a systematic way, in comparison to multiple state-of-the-art exact solvers for different problem classes.
    Our experiments show that despite being in its development stage, DA can provide high-quality solutions quickly, which, combined with the fact that the next generation of DA will allow up to one million binary variables \citep{da1MWeb}, indicates DA has the potential to become a successful technology in tackling large-scale COPs. 
    \item We find that DA outperforms the state of the art typically in large instances and is usually comparable otherwise, and illustrate that its performance is robust to the increase in the number of decision variables or constraints, as opposed to established solvers whose performance is known to be significantly sensitive to both.
    %
    \item We introduce two new applications for DA, namely QCPP and Sel-Col,
    and demonstrate that DA is particularly successful in yielding better solutions for instances that are challenging for the state-of-the-art solvers. 
    Moreover, we apply a size reduction heuristic for Sel-Col that 
    makes a significantly larger number of instances eligible for DA. 
    \item We analyze the effect of variable ordering and penalty coefficient values on the performance of DA, and show that they might need to be carefully selected in certain cases to increase the obtained solution quality.
    %
\end{enumerate}

The remainder of this paper is organized as follows. In Section \ref{Operating_principle}, we overview our solution approaches and the operating principles of state-of-the-art exact solvers and DA. 
In Section \ref{performance_analysis}, we analyze the performance of the solvers and DA on the four problem classes we consider. 
Finally, we provide conclusions and future research directions in Section \ref{conclusion}.

%% file: OperatingPrinciples.tex
\section{Solution approaches}
\label{Operating_principle}
In this section, we briefly explain the basic modeling and algorithmic principles behind three general methods to solve a COP, namely (i) solving an IP formulation of the problem via an exact solver, (ii) solving a QUBO formulation of the problem via an exact solver, and (iii) solving the QUBO formulation via DA.

\subsection{IP approach}
\label{approach:IP}
In general, any COP can be formulated as a binary linear program:
\begin{subequations}
\label{comb_model}
\begin{alignat}{2}
    \min \quad & c^\top x \label{m:COP:obj} \\
    \text{s.t.} \quad &  Ax =  b  && 	\label{m:COP:consts} \\
    & x \in \{ 0,1\}^n \label{m:COP:vardomain}
\end{alignat}
\end{subequations}
where $ c \in \R^{n}$, $A \in \R^{m \times n}$ and $b \in \R^{m}$. In other words, the goal is to find a solution $x \in \R^n$ which minimizes the linear objective function \eqref{m:COP:obj}, and satisfies all of the linear constraints \eqref{m:COP:consts} as well as variable domain restrictions \eqref{m:COP:vardomain}. 

State-of-the-art solvers such as GUROBI and CPLEX employ a linear programming-based branch-and-bound (BB) algorithm to solve binary linear programs. 
This algorithm conducts the search of an optimal solution in a structured way, specifically in the form of a search tree. 
At the root node of the tree, it solves the so-called linear programming relaxation of the model \eqref{comb_model} obtained by relaxing all integrality conditions on the decision variables, \eqref{m:COP:vardomain}. 
This yields a lower bound on the optimal value of the COP. 
If the obtained solution at the root node is not binary, the algorithm splits the search space (i.e., the relaxed feasible region) into smaller subregions based on branching rules, each excluding the infeasible relaxation solution. 
Creating one node per subregion in the search tree, the algorithm continues the search by repeating this procedure. 
Each node in the search tree represents the linear programming relaxation of model \eqref{comb_model} with additional constraints defining the subregion specific to that node.
If the solution at a node is binary, then it is a feasible solution for the original binary linear program, the associated objective value constitutes an upper bound for the optimal objective value, and there is no need to further branch on that node.
Otherwise, i.e., if there is a non-binary value in the solution, the algorithm decides on whether to further branch on that node by comparing the objective value of the node to the best upper bound found so far.
If the objective value is at least as large as the best upper bound, the algorithm discards that node because it cannot produce a better solution than the best one found so far.
In this manner, the BB algorithm avoids unnecessary computations in the search process.
Choosing well-suited branching rules and effective bounding techniques are two major considerations in the use of BB algorithms, which significantly impact the size of the search tree, thus the solution time.  
In addition to employing multiple carefully designed branching and bounding rules, the state-of-the-art solvers also implement a variety of enhancement techniques to improve the BB performance, including presolve methods to reduce the size of the models, cutting planes to strengthen the relaxations, and heuristics to discover good quality solutions during the search. 
For more information about recent developments on the BB algorithm, we refer the reader to \citep{morrison2016branch}.

It is also possible to incorporate a quadratic component in the objective function for certain problems, in which case model \eqref{comb_model} can be generalized into the following:
\begin{subequations}
\label{comb_model_quad}
\begin{alignat}{2}
    \min \quad & c^\top x + x^\top Q x \label{m:COP_quad:obj} \\
    \text{s.t.} \quad &  Ax =  b  && 	\label{m:COP_quad:consts} \\
    & x \in \{ 0,1\}^n, \label{m:COP_quad:vardomain}
\end{alignat}
\end{subequations}
\noindent where $Q \in \R^{n \times n}$ and symmetric.
Whether the objective function is linear or quadratic, we refer to models as in \eqref{comb_model} and \eqref{comb_model_quad} as IP models in this study, and call the approach of handling such models as the IP approach.

In order to solve binary quadratic problems, exact and heuristic methods exist in the literature. 
Initially proposed exact methods are (general) BB-based algorithms. 
In this regard, effective methods to make a quadratic model compatible with the BB framework and useful bounding techniques have been developed. 
Linearization methods aim to introduce a full linear representation of the quadratic model and tackle the problem using exact IP solvers and algorithms. 
This method requires the addition of new decision variables and constraints to the model, which in most cases increases the size of the problem significantly. 
On the other hand, outer approximation techniques try to approximate the objective function using supporting hyperplanes. 
They have been shown to be efficient for certain quadratic optimization problems, when supporting hyperplanes are easy to generate.
Exact solvers handle such models typically by benefiting from linearization and outer approximation techniques as well as other ideas used in treating nonlinear programs \citep{furini2013extended,gurobi2020nonconvex}.


There exist advanced solvers that can solve certain types of quadratic models exactly.
A mixed-integer quadratic programming (MIQP) model is one that can contain both continuous and integer variables, has linear constraints, and the objective function is allowed to contain quadratic terms. 
If the constraints are also allowed to contain quadratic terms, then the model is classified as a mixed-integer quadratically constrained programming (MIQCP) model.
Solvers that can handle MIQP and MIQCP models include CPLEX, GUROBI, SCIP, GLOMIQO, BARON, and MOSEK, with the last one being eligible except for the cases where quadratic objective/constraints and conic constraints are present.

We note that CPLEX and GUROBI are the state-of-the-art commercial solvers, most notably for mixed-integer programming, which are capable of efficiently solving certain types of quadratic problems optimally as well. 
In addition to these two powerful solvers, SCIP is a widely known open mixed-integer programming and mixed-integer nonlinear programming solver. 
In this study, we use these three solvers to compare against DA.

\subsection{QUBO approach}
\label{approach:QUBO}
We can reformulate and solve a COP as a QUBO model as well. For this purpose, we dualize the equality constraints, i.e., add them to the objective function along with a penalty coefficient. 
QUBO reformulation of model \eqref{comb_model_quad} is as follows:
\begin{subequations}
\label{QUBO_model}
\begin{alignat}{2}
    \min \quad & c^\top x  +  x^\top Q x  +  \penaltyCoef (Ax  -  b)^\top(Ax  -  b) \label{m:QUBO:obj} \\
    & x \in \{0,1\}^n, \label{m:QUBO:vardomain}
\end{alignat}
\end{subequations}
\noindent where $\penaltyCoef  > 0$ is a large number representing the unit penalty cost for constraint violation. 
We note that model \eqref{comb_model_quad} being a generalization of model \eqref{comb_model}, the QUBO reformulation \eqref{QUBO_model} is a valid one for binary linear model \eqref{comb_model} as well (by simply ignoring the $Q$-term in the objective function).
Since the constraints in model \eqref{comb_model_quad} are equalities, the dualized terms penalize any violation of an equality constraint in either direction in the objective function of model \eqref{QUBO_model}. 
If the penalty coefficient $\penaltyCoef$ is chosen ``sufficiently large", the model \eqref{QUBO_model} is an \emph{exact reformulation} of \eqref{comb_model_quad}, in the sense that any feasible solution to \eqref{comb_model_quad} yields a smaller objective value in \eqref{m:QUBO:obj} than any binary infeasible solution due to a large penalty term. 
Hence, if \eqref{comb_model_quad} has an optimal solution, then its optimal solution set and optimal objective value coincide with the ones of \eqref{QUBO_model}.

There exist exact and heuristic methods to solve QUBO problems, a comprehensive survey of which can be found in \citep{kochenberger2014unconstrained}.
Exact solvers for these problems typically use the techniques mentioned in Section \ref{approach:IP}.
QUBO problems are hard both in theory and practice, so that even very efficient exact algorithms can handle up to 500 variables \citep{mauri2012improving}.


\subsection{DA approach}
\label{approach:DA}
As an alternative to formulating IP models and solving them with state-of-the-art solvers,
some heuristic-based solution technologies that can handle QUBO models have been recently introduced. 
DA is one such emerging technology, which is a quantum-inspired hardware architecture and has an algorithm based on simulated annealing. DA has been shown to succeed in yielding high-quality (not necessarily optimal) solutions in a short amount of computation time for a variety of problems, including the routing and wavelength assignment problem \citep{cseker2020routing}, minimum vertex cover problem \citep{javad2019digitally}, maximum clique problem \citep{naghsh2019digitally}, and many more.

Simulated annealing is a local search method that proceeds by generating a candidate solution at each iteration and deciding whether to accept that solution.
If the proposed solution yields a better objective value than the current one, it is always accepted; otherwise, the acceptance decision is made based on the difference of the objective value of the candidate solution from that of the current, and a {\it temperature} parameter.
The temperature regulates the tendency of the algorithm to accept non-improving candidate solutions and is typically non-increasing as the search progresses.
Higher temperatures more likely let the algorithm explore a larger region of the objective function by increasing the acceptance probability of non-improving solutions and thereby help escape from local optima,
while the search intensifies around a narrower area as the temperature decreases.

The algorithm of DA is based on simulated annealing, but it incorporates multiple additional features by utilizing its computational architecture, such as parallel evaluation of all neighbouring solutions as well as some techniques to avoid local optima.
We next describe the operating principles of DA, mostly based on \citep{DAU2}.

DA tries to minimize a quadratic function of binary variables, a so-called energy function $\E$, which corresponds to the objective function of a QUBO model, as in \eqref{QUBO_model}. 
If we let $\hat{x} = (\hat{x}_1, \hat{x}_2,...,\hat{x}_n)$ be the current solution at an iteration, a candidate solution is constructed by flipping the value of one variable from $\hat{x}_j$ to $1-\hat{x}_j$.
At each iteration, DA considers all possible candidate solutions by flipping the value of each variable separately, and calculates the changes in the energy function value in parallel, i.e., it computes $\Delta \E_{j}$ values for all $j = 1,...,n$ simultaneously.
The acceptance condition it uses is the so-called Metropolis criterion, which sets the probability of accepting the flip as
\begin{equation}
    Pr(\Delta \E_{j}) = \min \left(1, \exp(-\beta  \cdot   \Delta \E_j)\right),
\label{prob}
\end{equation}
where $\beta$ is the inverse of the temperature parameter. 
For each $j$, a binary flag becomes ``1" with probability \eqref{prob}, meaning that the flip can be accepted. 
Among all variables with flag ``1", DA chooses one and proceeds to the next state by flipping the value of the associated variable.

Throughout the search procedure, DA might reach to a local minimum state, where no candidate variable to flip can be found. 
In such cases, a mechanism to escape local optima is facilitated by adding a positive offset $\E_{\text{off}}$ to the energy function, equivalent to multiplying the acceptance probabilities with $\exp(\beta \ \E_{\text{off}})$.
So, when there is no variable with flag ``1", the offset generator dynamically adds a constant to $\E_{\text{off}}$; otherwise, the offset value is set to zero.

In addition to the simultaneous evaluation of all candidate solutions and the dynamic offset mechanism, which mainly differentiate DA's algorithm from simulated annealing, 
DA has the {\it parallel tempering} option, also referred to as the {\it replica exchange} method, where multiple independent search processes with different initial temperatures function in parallel and exchange solutions with each other. 
In this study, we refer this mode of operation as the {\it parallel mode} of DA, and the regular one as the {\it normal mode}.
We utilize the second generation of DA in our experiments, which can handle up to 8192 binary variables. 

%% file: PerformanceAnalysis.tex
\section{Performance analysis}\label{performance_analysis}


In this section, we provide the results of our extensive experimental study that we conducted for comparing the performance of DA \citep{da2Web,DAU2} and three state-of-the-art solvers, namely GUROBI \citep{optimization2018}, CPLEX \citep{cplex2009v12}, and SCIP \citep{Achterberg2009}, on different problem classes and their benchmark instances from the literature and/or newly generated ones. 
Since DA inherently uses a heuristic (rather than an exact) algorithm, comparing the solutions produced by DA to those of solvers, which can yield provably optimal solutions, is an integral part of this study to gain insight into the quality of DA's solutions.

We conduct our experiments with DA using the Digital Annealer environment prepared exclusively for research at the University of Toronto.
We run our experiments with GUROBI, CPLEX and SCIP on a Mac computer with 3 GHz Intel Core i5 CPU and 16 GB memory, using up to four threads.
We use the latest release of each solver we consider. 
In particular, we use GUROBI 9.0.3, CPLEX 12.10, SCIP 7.0.1, and the second generation of DA. 
We experiment with both the parallel and normal modes of DA, using the default hyperparameters for the former, and with a tuned set of run parameters for the latter, which we fix as a result of some preliminary tests.


In the assessment of experimental outcomes, we use various measures for a clear and explanatory illustration, as the nature of the problem and obtained results necessitate. 
One widely used performance measure for optimization problems is the percentage optimality gap, defined as 
\begin{align}
    \text{\% solver gap} =  100 \cdot \frac{\text{UB}-\text{LB}}{|\text{UB}|}, \label{def:solvergap}
\end{align}
where UB and LB respectively denote the lower and upper bound values on the optimal objective value. 
This value measures how far the objective value of the best feasible solution found can be from the optimal value, and indicates how successful the solver is in proving the optimality of a solution.
In order to gain better insight into the closeness of an objective value to the optimal, rather than the success of the solver in proving the goodness of a solution, we can calculate percentage gap values using the best available lower bound, as follows:
\begin{align}
    \text{\% gap} =  100 \cdot \frac{\text{UB}-\text{LB}_{\text{best}}}{|\text{UB}|}, \label{def:gap}
\end{align}
\noindent where $\text{LB}_{\text{best}}$ is the maximum lower bound value among a set of available values (e.g., those from different solvers under different time limits).
The percentage optimality gap definition as in \eqref{def:solvergap} is the widely used one, and we refer to it as {\it \% solver gap} in this paper to differentiate from the one that is calculated with the best available lower bound value, provided in \eqref{def:gap}, which we simply call {\it \% gap}.

Despite the fact that percentage optimality gap is a widely used performance measure, it is not one that we predominantly utilize in this study, because even though the three solvers other than DA declare themselves to be exact in solving binary quadratic problems, we observed during our computational study that some occasionally yield invalid lower bounds in these cases (assuming problems are in minimization form), i.e., lower bounds that exceed the upper bounds.
In order to guarantee the validity of the utilized performance measures, we do not base our analysis on measures involving lower bounds for quadratic cases.
Instead, we work with normalized difference of upper bound values, which we define as
\begin{align}
    \text{Normalized difference} =  100 \cdot \frac{\text{UB}-\text{UB}_\text{ref}}{|\text{UB}_\text{ref}|}, \label{def:normdiff}
\end{align}
where $\text{UB}_\text{ref}$ is a reference upper bound value we normalize with respect to.
We typically set the $\text{UB}_\text{ref}$ values as the upper bounds from a promising solver and time limit combination.
Since all of the problems we consider in this study are in minimization form, upper bound values simply correspond to the objective values of feasible solutions, and we use the two terms interchangeably.
We introduce other performance measures in the upcoming sections, as needed.

We note that one cannot directly impose a time limit on DA, but can set the run parameters, e.g., the number of iterations, in such a way that a certain level of solution time is achieved.
Once the run parameters are fixed, the solution times stay roughly the same throughout all instances for a given problem type.
Also, the solution times of DA do not typically exceed few minutes, even with very high numbers of iterations.
Therefore, we mostly make our comparisons on the basis of up to two-minute run times, specific values of which we vary for different problem classes as needed.
When we discuss the results of our timed experiments in the sections that follow, we simply mean that we have fixed the run parameters of DA to achieve a certain solution time, and explicitly imposed time limits for the other solvers.

In the sequel, we first focus on pure QUBO problems (Section \ref{subsec:pureQUBO}), i.e., problems that naturally appear in QUBO form.
We then 
extend our study to constrained cases and consider the quadratic assignment (Section \ref{subsec:QAP}), quadratic cycle partition (Section \ref{subsec:QCPP}), and selective graph coloring problems (Section \ref{subsec:selcol}).


\input{PureQUBO}

\input{QAP}

\input{QCPP}

\input{SelCol}

%% file: PureQUBO.tex
\subsection{Pure QUBO}
\label{subsec:pureQUBO}

In this section, we concentrate on pure QUBO problems; that is, problems that are formulated as a QUBO model in the literature, equivalent to model \eqref{QUBO_model} with $\penaltyCoef=0$, rather than being recast into this structure through some transformation. 

Various heuristic and exact methods have been developed and tested for QUBO models. 
Heuristic methods include tabu search methods \citep{glover1998adaptive,glover1999tabu,glover2010diversification}, local search methods \citep{boros2007local}, and simulated annealing \citep{beasley1998heuristic,hasan2000comparison}.
Existing exact methods include
variants of the branch-and-bound algorithm \citep{huang2006lower,helmberg1998solving, billionnet1994minimization}, Lagrangian decomposition \citep{mauri2011lagrangean}, and linearizations of pure QUBO problems \citep{gueye2009linearization}.
More information on QUBO applications and solution approaches can be found in the survey by \cite{kochenberger2014unconstrained}.

\subsubsection{Instances}

We collected pure QUBO instances of varying sizes from well-known libraries; namely, OR-Library\footnote{
\href{http://people.brunel.ac.uk/~mastjjb/jeb/orlib/bqpinfo.html}{http://people.brunel.ac.uk/~mastjjb/jeb/orlib/bqpinfo.html}
} \citep{beasley1990or}, Biq-Mac\footnote{
\href{http://biqmac.uni-klu.ac.at/biqmaclib.html}{http://biqmac.uni-klu.ac.at/biqmaclib.html}
} \citep{wiegele2007biq}, and QPLIB\footnote{
\href{http://qplib.zib.de}{http://qplib.zib.de}
} \citep{furini2018qplib}. 
From OR-Library, we use test instances in \citep{beasley1998heuristic}, which have 20 to 2500 variables. 
Biq-Mac is a library that contains medium-size QUBO problem instances with 20 to 500 variables, some of them being the same as those offered by OR-library. Finally, QPLIB is a comprehensive library for quadratic programming that involves 23 pure QUBO problems with 120 to 1225 variables.
Using these three libraries, we obtain 208 pure QUBO instances in total which contain 20 to 2500 variables. 
All of these instances are eligible for DA. 
In order to additionally test some larger instances that DA can handle, we generated 55 new ones having 3000 to 8000 variables, following the parameters used in generating \cite{beasley1990or} type instances available in the OR-Library and also Biq-Mac library.
Pure QUBO test bed information is summarized in 
Table \ref{tab:pureQUBO_instance_info}. 
In our experiments, for 3 of the 23 instances from QPLIB, instance parameters turned out to be beyond the acceptable ranges for DA.
Therefore, we present the results for the remaining 260 instances in total.

\input{tables/pureQUBO_instance_info}

\subsubsection{Experimental results}

\begin{figure}[ht]
    \centering
	\begin{subfigure}{0.49\textwidth}
        \centering
        \includegraphics[scale=0.62]{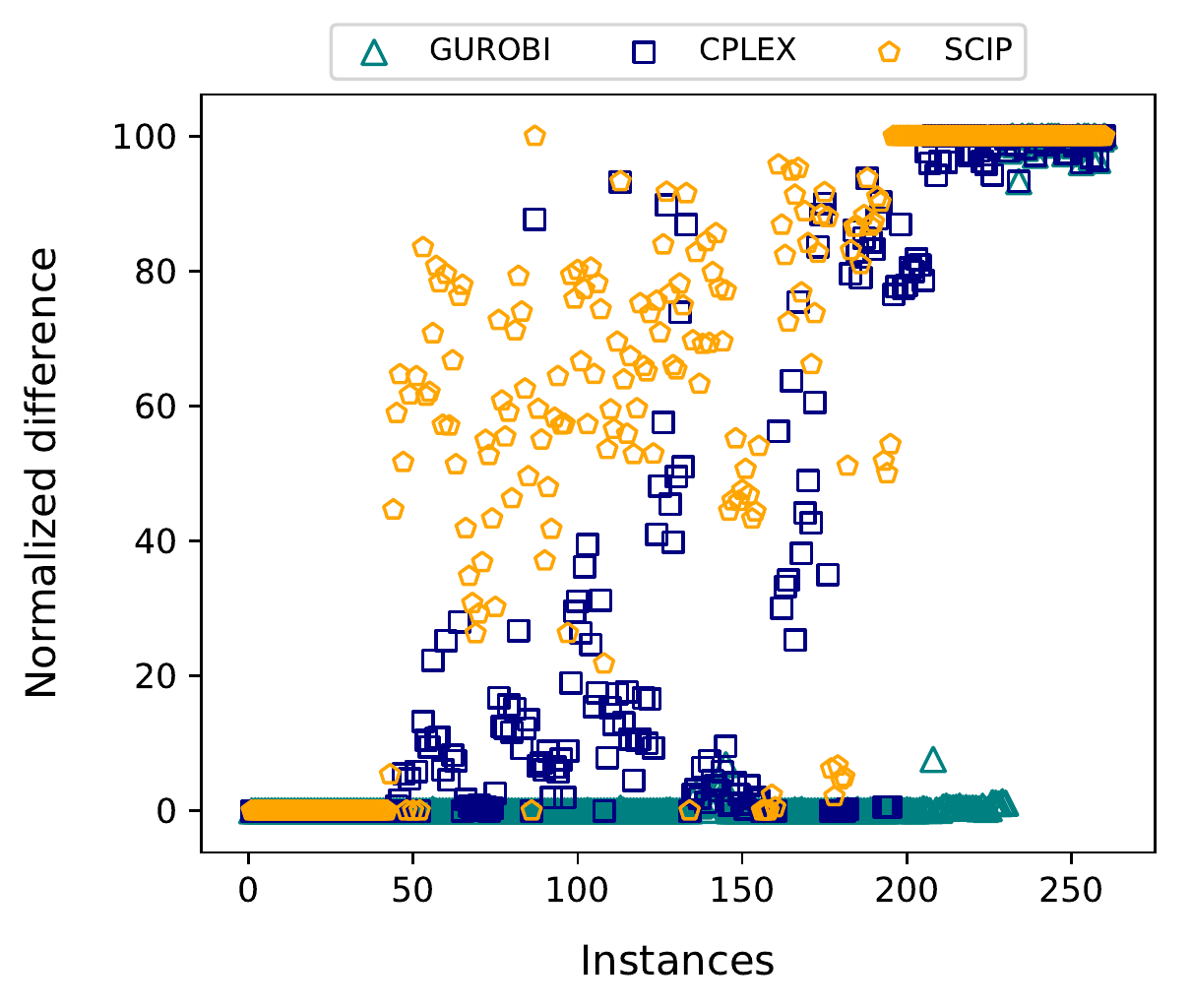}
        \caption{60-second time limit}
        \label{fig:pureQUBO_exact_solvers_60s}
	\end{subfigure}
	\begin{subfigure}{0.49\textwidth}
		\centering
	    \includegraphics[scale=0.62]{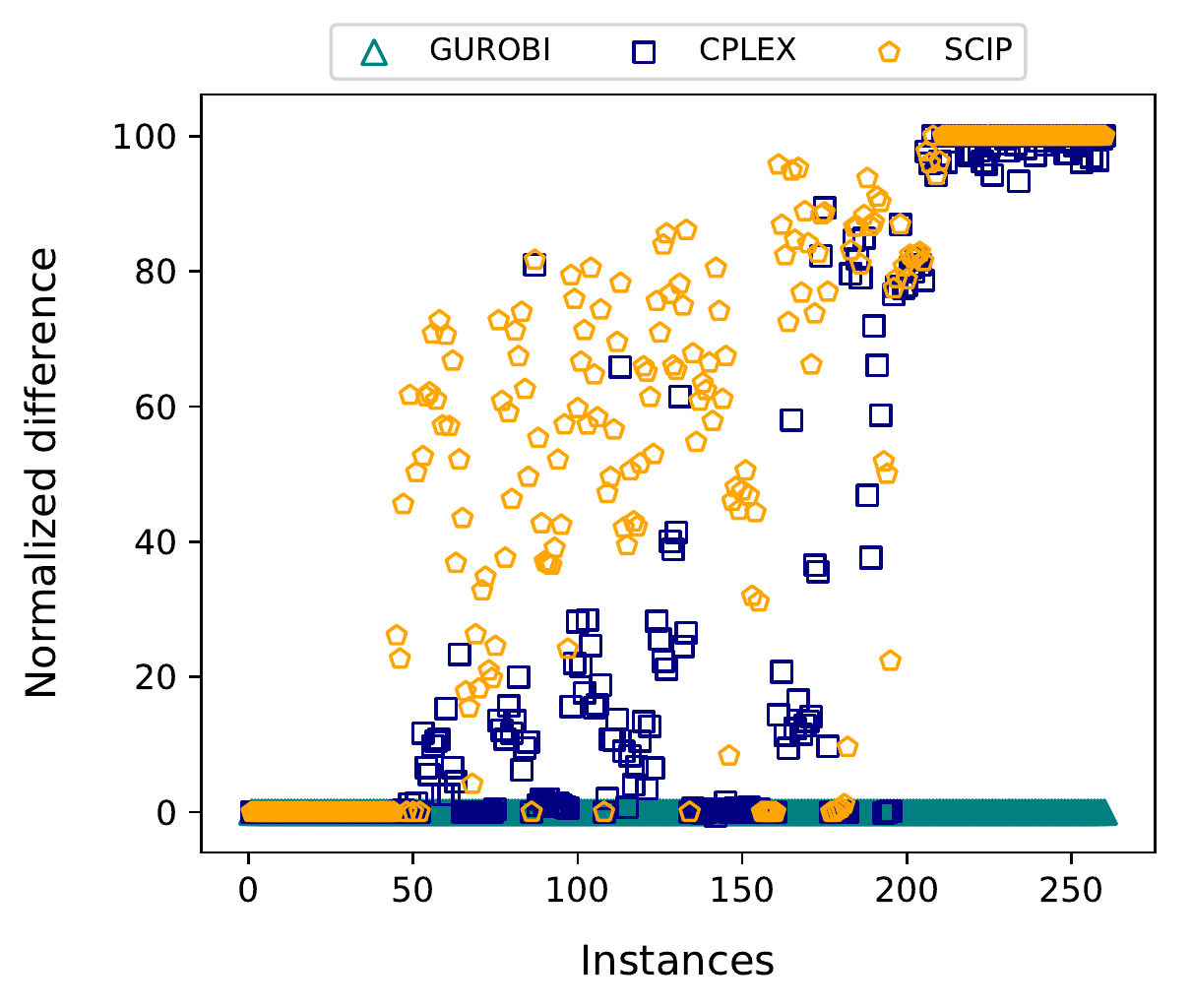}
		\caption{10-minute time limit}
		\label{fig:pureQUBO_exact_solvers_10min}
	\end{subfigure}
	\caption{Performance analysis of exact solvers for pure QUBO instances, with GUROBI 10-minute results taken as reference.}
	\label{fig:pureQUBO_exact_solvers}
\end{figure}

{\partitlestyle Comparison of exact solvers.} We begin with the comparison of the three exact solvers; GUROBI, CPLEX, and SCIP. 
Figure \ref{fig:pureQUBO_exact_solvers} displays a comparison of the performances of the three solvers, under 60-second (Figure \ref{fig:pureQUBO_exact_solvers_60s}) and 10-minute (Figure \ref{fig:pureQUBO_exact_solvers_10min}) time limits. 
The $x$-axes represent the instances sorted based on the number of variables they contain, which will be so in all the figures that we present, while the $y$-axes show the normalized differences of the upper bound values of the solvers, with the upper bound values from GUROBI under 10-minute time limit taken as reference, unless otherwise stated. 

Under 60-second time limit, GUROBI outperforms both CPLEX and SCIP for small and moderate-size instances. 
As the instances become larger, the upper bounds provided by GUROBI and CPLEX become closer. 
When we increase the time limit to 10 minutes, the solution quality of GUROBI considerably improves especially for large instances, while CPLEX and SCIP produce almost the same results as before.

\begin{figure}[ht]
    \centering
	\begin{subfigure}{0.49\textwidth}
	    \centering
        \includegraphics[scale=0.62]{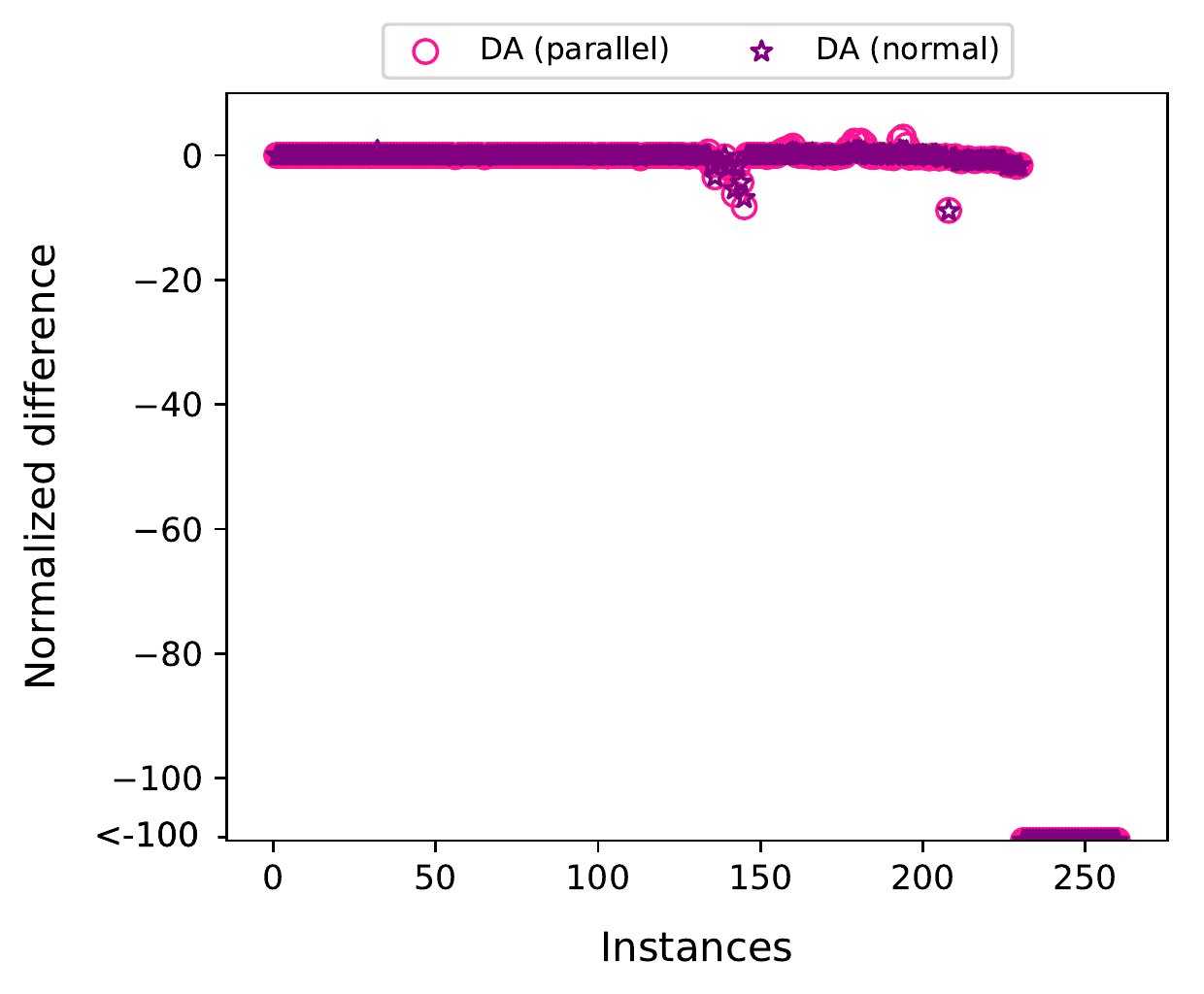}
        \caption{GUROBI 60-second as the reference}
        \label{fig:pureQUBO_da_vs_da_60s}
	\end{subfigure}
	\begin{subfigure}{0.49\textwidth}
        \centering
	    \includegraphics[scale=0.62]{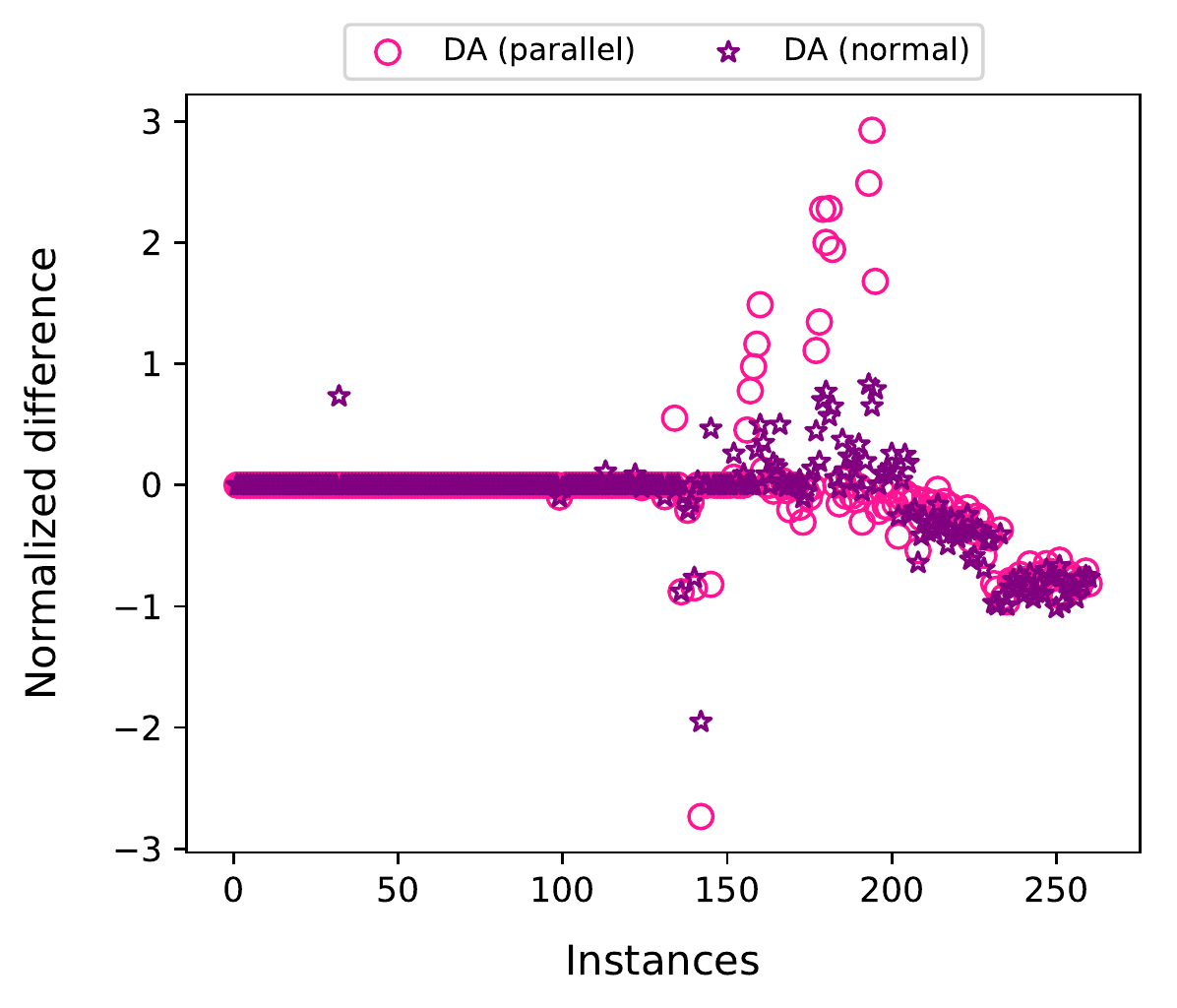}
		\caption{GUROBI 10-minute as the reference}
		\label{fig:pureQUBO_da_vs_da_10_min}
	\end{subfigure}
	\caption{Performance comparison of normal and parallel modes of DA for pure QUBO instances, under 60-second time limit.}
	\label{fig:pureQUBO_da_vs_da}
\end{figure}

{\partitlestyle Comparison of DA modes.} 
Next, we evaluate the performances of normal and parallel modes of DA.
When running the normal mode of DA, we used the best parameter combination we found after some preliminary testing, although we do not claim it to be the best possible of all, {\cfujitsu and for the parallel mode, we used the default hyperparameters, as mentioned at the beginning of Section \ref{performance_analysis}}.
Figure \ref{fig:pureQUBO_da_vs_da} shows the performance of DA in terms of normalized differences, with GUROBI's upper bound values from 60-second and 10-minute experiments being used as reference, provided in Figures \ref{fig:pureQUBO_da_vs_da_60s} and \ref{fig:pureQUBO_da_vs_da_10_min}, respectively. 


 

We observe that
the two modes of DA yield almost the same level of performance, with the normal mode being slightly better than parallel. 
Moreover, 
the average running time of the normal mode (0.40 second) is significantly less than that of parallel mode (67.03 seconds). 
Since these plots are constructed by taking the upper bound values from GUROBI as a reference, they also implicitly compare the performances of DA and GUROBI under 60-second and 10-minute time limits. 
We see from Figure \ref{fig:pureQUBO_da_vs_da_60s} that for small and moderate-size instances, both modes of DA and GUROBI provide similar upper bounds under 60-second time limit. 
As the instance size increases, the upper bounds provided by DA become significantly better than those of GUROBI. 
Though GUROBI is able to yield improved solutions when the time limit is increased to 10 minutes, for large instances it is still outperformed by DA.
So, we conclude that DA is particularly effective in tackling large instances where exact solvers start struggling.

{\partitlestyle Variable ordering analysis.} 
We now analyze the effect of input format on the performance of DA, namely the impact of different variable orderings, inspired by the observation that it can significantly affect the performance of established solvers in solving mixed-integer linear programs (see for instance \citep{achterberg2013mixed}). 
For this purpose, we take a sample of five pure QUBO problem instances with 1000 variables, and experiment with 100 different random permutations of variable indices, using both modes of DA. 
For normal mode, we test two different parameter combinations; an arbitrarily selected one that is not specifically tuned, and one that is the best we obtained as a result of some preliminary experiments.

Table \ref{tab:pureQUBO_var_ordering} summarizes the results of our tests with different variable orderings. 
For each one of the three alternative configurations, namely the parallel mode (``Parallel"), normal mode with an arbitrary parameter configuration (``Normal"), and normal mode with a tuned parameter combination (``Normal (tuned)"), we provide the average of the best found objective values for 100 random permutations of variable indices. 
In order to examine possible effects of the number of iterations, we conduct this experiment using two different settings for it; $10^3$ and $10^6$ iterations. 
In the last row of the table, we also report the average percentage difference (``Avg \% diff") of the objective values corresponding to the best and worst orderings, where the differences are normalized with respect to the worst ordering.

\input{tables/pureQUBO_var_ordering}

The results indicate that while the impact of variable ordering is insignificant when the number of iterations is large, it does make a difference when the number of iterations is low, as the average percentage differences reveal.
Additionally, we observe that when the parameters for normal mode are not tuned, there is room for improvement through different variable orderings. 
In particular, with a relatively small number of iterations, the difference between the best and worst variable ordering decreases from about 42\% to 25\% with the tuning of parameters, which is a significant improvement. 
So, we conclude that, when the number of iterations is high, DA becomes insensitive to the input format and demonstrates a more robust performance.

{\partitlestyle Overall performances.} 
Finally, we present a summary of our computational results on pure QUBO instances in Table \ref{tab:pureQUBO_summary}, which contains average normalized differences with 10-minute results of GUROBI taken as reference (``Norm diff"), average percentage gaps with respect to the best lower bound found (``\% gap"), average percentage gaps reported by solvers (``\% solver gap"), as well as average run times (``Time (sec)"). 
The \largeVal~signs in the table stand for values that are greater than $10^5$.
The results indicate that DA both in normal and parallel modes outperforms the other solvers both in terms of solution quality and run times. 

\input{tables/pureQUBO_summary}

%% file: tables/pureQUBO_instance_info.tex
\begin{table}[htbp]
\centering
	\caption{Information about pure QUBO instances.}
	\label{tab:pureQUBO_instance_info}
	\def\arraystretch{1}
	 \scalebox{0.86}{
    \begin{tabular}{lrr}
    	\toprule
    	& \parbox[t]{2.5cm}{\centering \# instances} & \parbox[t]{2cm}{\centering  \# variables}\\
    	\midrule
    	{OR-Library} & 105\phantom{ooo} & 20--2500\\
    	{Biq-Mac} & 80\phantom{ooo} & 100--500\\
    	{QPLIB} & 23\phantom{ooo} & 120--1225 \\
    	New & 55\phantom{ooo} & 3000--8000 \\ 
        \midrule
        {Total} & 263\phantom{ooo} & 20--8000\\
    	\bottomrule
    \end{tabular}%
    }
\end{table}

%% file: tables/pureQUBO_var_ordering.tex

\begin{table}[ht]
	\centering
	\caption{Results of experiments on the impact of variable ordering for DA.}
	\label{tab:pureQUBO_var_ordering}
	 \def\arraystretch{1}
	  \scalebox{0.83}{
	\begin{tabular}{c S[table-format=5.1] S[table-format=4.1] S[table-format=5.1] S[table-format=6.1] S[table-format=6.1] S[table-format=6.1]}
		\toprule	
		 & &\multicolumn{1}{c}{$10^3$ iterations} & & &\multicolumn{1}{c}{$10^6$ iterations} & \\
		 \cmidrule(lr){2-4}\cmidrule(lr){5-7}
		\parbox[t]{2.25cm}{\centering Instance} 
		& \parbox[t]{1.9cm}{\centering Normal} 
		& \parbox[t]{1.9cm}{\centering Normal} 
		& \parbox[t]{1.9cm}{\centering Parallel} 
		& \parbox[t]{1.9cm}{\centering Normal} 
		& \parbox[t]{1.9cm}{\centering Normal} 
		& \parbox[t]{1.9cm}{\centering Parallel} \\
		& \multicolumn{1}{c}{(tuned)} & \multicolumn{1}{c}{} & \multicolumn{1}{c}{} & \multicolumn{1}{c}{(tuned)} & \multicolumn{1}{c}{} & \multicolumn{1}{c}{} \\
		\midrule
		bqp\_1000\_1 & -17187.4 & -8413.3  & -16313.2 & -185495.3 & -185791.0  & -185791.0 \\[0.05cm]
		bqp\_1000\_2 & -17258.5 & -8421.9  & -16176.1 & -176854.3 & -177526.0  & -177526.0 \\[0.05cm]
		bqp\_1000\_3 & -17358.8 & -8580.2  & -16499.1 & -185094.3 & -185575.0  & -185575.0 \\[0.05cm]
		bqp\_1000\_4 & -17417.0 & -8614.5 & -16449.0 & -184899.9 & -185253.0  & -185253.0 \\[0.05cm]
		bqp\_1000\_5 & -17156.9 & -8334.3  & -16178.4 & -175915.8 & -176566.0  & -176566.0 \\[0.05cm]
		\midrule
		 Avg \% diff &  25.0 &  41.9 & 26.3 & 0.5 &  0.0 & 0.0 \\[0.05cm]
		\bottomrule
	\end{tabular}%
	}
\end{table}%

%% file: tables/pureQUBO_summary.tex
\begin{table}[htbp]
  \centering
  \caption{Average performance of DA, GUROBI, CPLEX and SCIP for pure QUBO instances.}
  \def\arraystretch{1}
  \scalebox{0.72}{
    \begin{tabular}{crrrrrrrrrrrr}
    \toprule
      &       & \multicolumn{4}{c}{60-sec limit} & \multicolumn{4}{c}{120-sec limit} & \multicolumn{3}{c}{10-min limit} \\
      \cmidrule(lr){3-6}\cmidrule(lr){7-10}\cmidrule(lr){11-13}
      & \parbox[t]{1.3cm}{\centering DA (n)} 
      & \parbox[t]{1.3cm}{\centering DA (p)} 
      & \parbox[t]{1.4cm}{\centering GUROBI} 
      & \parbox[t]{1cm}{\centering CPLEX}
      & \parbox[t]{1cm}{\centering SCIP}  
      & \parbox[t]{1.3cm}{\centering DA (p)}  
      & \parbox[t]{1.4cm}{\centering GUROBI} 
      & \parbox[t]{1cm}{\centering CPLEX} 
      & \parbox[t]{1cm}{\centering SCIP} 
      & \parbox[t]{1.4cm}{\centering GUROBI} 
      & \parbox[t]{1cm}{\centering CPLEX} 
      & \parbox[t]{1cm}{\centering SCIP} \\
    \midrule
    Norm diff & -0.1 & -0.1 & 11.6 & 36.9 & 60.2 & -0.2 & 9.7 & 34.3 & 58.0  & 0.0 & 32.3 & 55.1 \\[0.05cm]
    \% gap & 9.2 & 9.2 &  \largeVal &  \largeVal &  \largeVal & 9.1 & \largeVal & \largeVal & \largeVal & 9.3 & \largeVal & \largeVal \\[0.05cm]
    \% solver gap &  {--} &  {--} &  \largeVal &  \largeVal &  \largeVal &  {--} &  \largeVal &  \largeVal &  \largeVal & 9.9 & \largeVal & \largeVal \\[0.05cm]
    Time (sec) & 0.4 & 67.0 & 48.5 & 48.1 & 50.1 & 112.6 & 96.2 & 96.2 &  98.5 & 472.4 & 487.8 & 480.0 \\[0.05cm]
    \bottomrule
    \end{tabular}%
    }
  \label{tab:pureQUBO_summary}%
\end{table}%

%% file: QAP.tex
\subsection{Quadratic assignment problem}
\label{subsec:QAP}

In this section, we focus on the well-known quadratic assignment problem, QAP.
It is typically formulated as a constrained model, but can easily be reformulated as a QUBO model as well, because the original constraints are equalities and all the decision variables are binary.

QAP is one of the most difficult COPs and has appeared in many practical applications, 
including economic problems, scheduling, hospital planning, placement of electronic components, and many more. 
Different mathematical formulations (e.g., integer and mixed-integer linear, permutation-based, graph-based formulations), heuristics (e.g., local search, tabu search, simulated annealing, genetic algorithms) and exact solution approaches (e.g., branch-and-bound, cutting plane algorithms, dynamic programming) have been proposed for it.
For more information on QAP, we refer the reader to the survey by \cite{loiola2007survey}.

In what follows, we consider the IP formulation initially proposed by \cite{koopmans1957assignment}, reformulate it as a QUBO model, and present the results of our experiments using these models with DA and the three solvers.

\subsubsection{Formulations}
We are given a set of $n$ facilities and $n$ candidate locations. For each pair of locations $ \{k,\ell\} $, a distance parameter $d_{k\ell}$, and for each pair of facilities $ \{i,j\} $, a flow parameter $f_{ij}$ is defined. 
The goal is to assign each facility to a distinct location by minimizing the sum of the product of distances with the corresponding flows. 
To this end, we define a set of binary decision variables: 
\begin{equation*}
    x_{ik} = 
    \begin{cases}
      1, & \text{if facility $i$ is assigned to location $k$} \\
      0, & \text{otherwise}.
    \end{cases}
\end{equation*}

Using the parameters and variables as defined above, an IP formulation for QAP can be written as follows:
\begin{subequations}
\begin{alignat}{2}
    \min \quad & \sum_{i=1}^{n} \sum_{j=1}^{n} \sum_{k=1}^{n} \sum_{\ell=1}^{n} f_{ij} \ d_{k\ell} \ x_{ik} \ x_{j\ell} &&  \label{eq:QAP_obj} \\[0.12cm]
    \text{s.t.} & \sum_{i=1}^{n} x_{ik} = 1 && k \in \{1,\hdots,n\}  \label{eq:QAP_c1} \\
    & \sum_{k=1}^{n} x_{ik} = 1    && i \in \{1,\hdots,n\}  \label{eq:QAP_c2} \\
    & x_{ik} \in \{0,1\}   && i \in \{1,\hdots,n\}, \ k \in \{1,\hdots,n\} \label{eq:QAP_c3}	
\end{alignat}
\end{subequations}
The quadratic objective function \eqref{eq:QAP_obj} minimizes the total cost of the assignment. Constraints \eqref{eq:QAP_c1} ensure that each location is assigned exactly one facility, and \eqref{eq:QAP_c2} guarantee that each facility is assigned to a single location. Finally, constraints \eqref{eq:QAP_c3} enforce the decision variables to be binary.

Dualizing the equality (assignment) constraints, we obtain the following QUBO formulation of QAP:
\begin{subequations}
\begin{align}
\min\quad &{  \sum_{i=1}^{n} \sum_{j=1}^{n} \sum_{k=1}^{n} \sum_{\ell=1}^{n}  f_{ij}  d_{k\ell}  x_{ik}  x_{j\ell} }  +  \penaltyCoef \sum_{k=1}^{n}\left(\sum_{i=1}^{n} x_{ik} - 1 \right)^2  +  \penaltyCoef \sum_{i=1}^{n}\left(\sum_{k=1}^{n} x_{ik} - 1 \right)^2  \label{eq:QAP_qubo_obj}\\[0.25cm]
\text{s.t.} \quad & x_{ik} \in \{0,1\} \qquad  i \in \{1,...,n\}, \ k \in \{1,...,n\} \label{eq:QAP_qubo_c}
\end{align}
\end{subequations}

\noindent where $\penaltyCoef > 0$ is the penalty coefficient. We note that, for a sufficiently large value of $\penaltyCoef$, the QUBO formulation for QAP provided in \eqref{eq:QAP_qubo_obj}--\eqref{eq:QAP_qubo_c} is exact; i.e., it produces provably optimal solutions for the original problem.

\subsubsection{Experimental results}


We collected 136 QAP instances from the well-known library named QAPLIB\footnote{
\href{http://qplib.zib.de}{http://qplib.zib.de}
} \citep{furini2018qplib}. 
The instances have 100 to 10000 variables. 
After filtering the instances that are size-wise eligible for DA, a total of 123 instances remained. 
We conducted our computational experiments on those instances, using the three state-of-the-art solvers to solve the associated IP and QUBO formulations, and DA for solving the QUBO formulations.
Since it takes less than two minutes for DA to solve each of these instances, we used 60- and 120-second time limits in our experiments to compare the performances of solvers on equal terms.
In addition, we conducted 10-minute experiments using GUROBI, CPLEX and SCIP in order to see how much their solution quality improves when a longer time is allowed.
We used a penalty coefficient ($\penaltyCoef$) value of 16,000 in our QUBO experiments, which, among the ones we tested, led to the best overall results for DA.
In particular, higher $\penaltyCoef$ values, combined with the instance coefficients in our test bed, went beyond the acceptable ranges for DA, where smaller values led to lower feasibility percentages.

{\partitlestyle Feasibility analysis.} 
Our first observation from our experiments is that, when solving the QUBO formulations, the solvers fail to yield feasible solutions for a significant portion of the QAP instances.
Therefore, we begin with analyzing how much the solvers succeed in yielding feasible solutions.
Table \ref{tab:QAP_feasibility} reports the feasibility percentage of the solutions obtained within 120-second time limit (``\% feas"), as well as the average percentage gap of the solvers on the instances for which DA fails to find any feasible solutions (``\% solver gap"). 

\input{tables/QAP_feasibilityTable}

From the first row of the table, we observe that GUROBI, CPLEX and SCIP are able to find feasible solutions to the IP formulations of all QAP instances within the time limit. 
However, as QUBO solvers, they fail to find any feasible solutions for a considerably large subset of instances. 
Among QUBO solvers, both normal and parallel modes of DA outperform the others in providing feasible solutions. 
Furthermore, the second row of the table indicates that although exact IP solvers can find feasible solutions for all of the instances, their average percentage gaps are too high in the instances for which DA fails to find any feasible solutions. 
Moreover, GUROBI, CPLEX and SCIP as QUBO solvers cannot find feasible solutions for all of those instances within the time limit.
These findings imply that those instances are too hard thus the solvers struggle in finding high-quality solutions.

In order to gain further insight into the feasibility percentage of the solutions DA yields, we chose a subset of instances and scaled down the coefficients of the distance and flow matrices by dividing them by a constant number. 
We selected this subset in such a way that it contains instances where DA performs well, moderately and poorly, so that the results obtained on this test set can be reliably generalized.
We have seen that with the new scaling, the feasibility percentage in DA's solutions on the selected set of instances increases from about 34\% to 75\%, implying that the performance of DA is sensitive to the magnitude of the input coefficients.

In the sequel, we present the experimental results obtained on the original QAP instances. 
Each figure shows the results from a list of common instances eligible for a particular comparison, with the $x$-axes representing instances sorted in ascending order with respect to the number of variables they contain.
We begin with comparing the modes of DA, and then solvers GUROBI, CPLEX, and SCIP to each other both as IP solvers and QUBO solvers. 

\begin{figure}[h]
    \centering
	\includegraphics[scale=0.62]{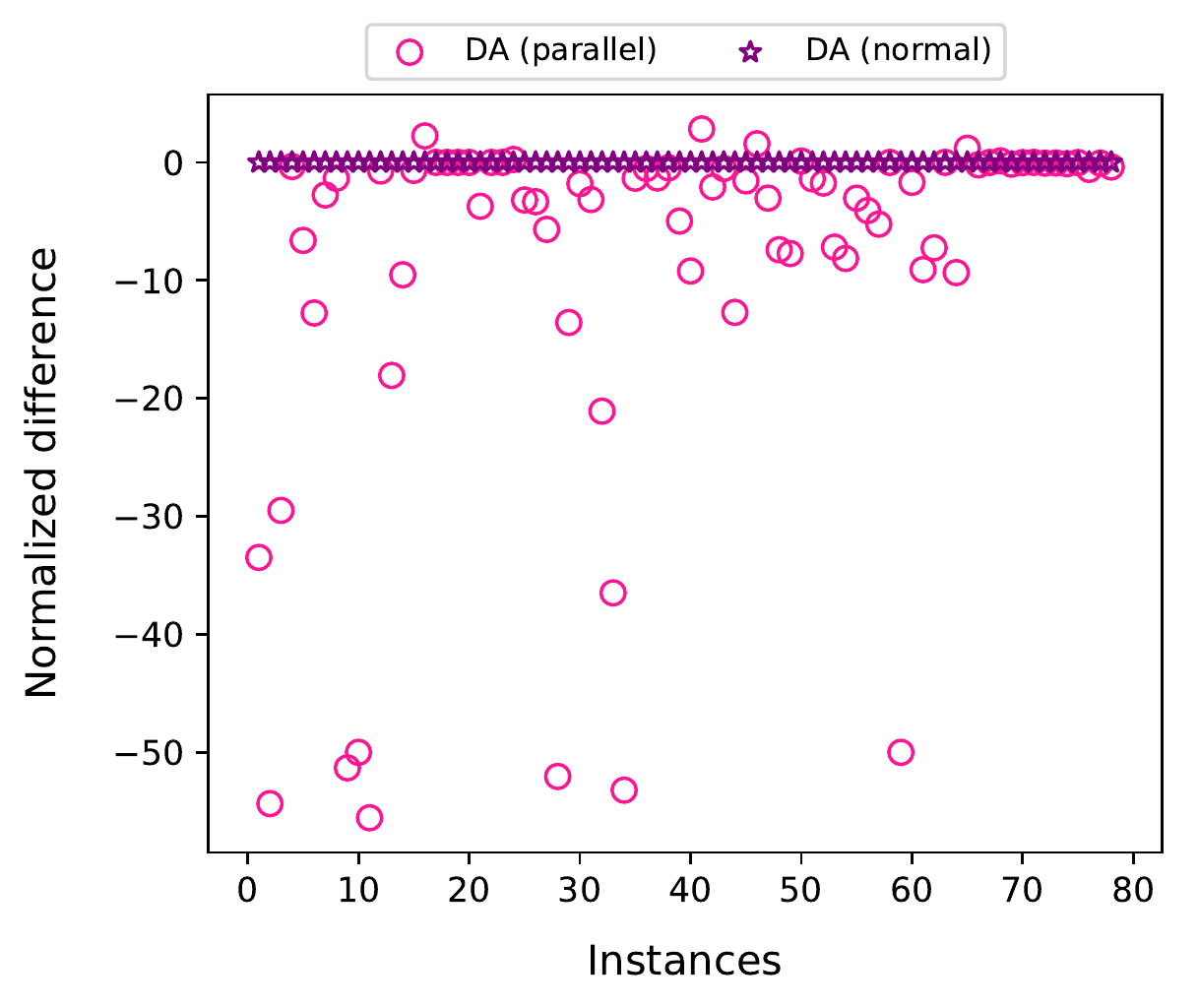}
	\caption{Performance comparison of DA modes for QAP 
	under 120-second time limit, with the results from the normal mode of DA taken as reference.}
	\label{subfig:QAP_da_vs_da}
\end{figure}

{\partitlestyle Comparison of DA modes.} 
Figure \ref{subfig:QAP_da_vs_da} illustrates the performance of DA in normal and parallel modes on the basis of normalized difference of the upper bound values.
Here, the objective values from the normal mode of DA is taken as reference in order to fix one set of data points to zero level and thereby facilitate the comparison. 
We observe that parallel mode of DA outperforms normal mode especially on small and medium-size instances, because the majority of the data points for the parallel mode of DA lies below the zero level. 

\begin{figure}[h]
    \centering
	\begin{subfigure}{0.49\textwidth}
		\centering
	    \includegraphics[scale=0.62]{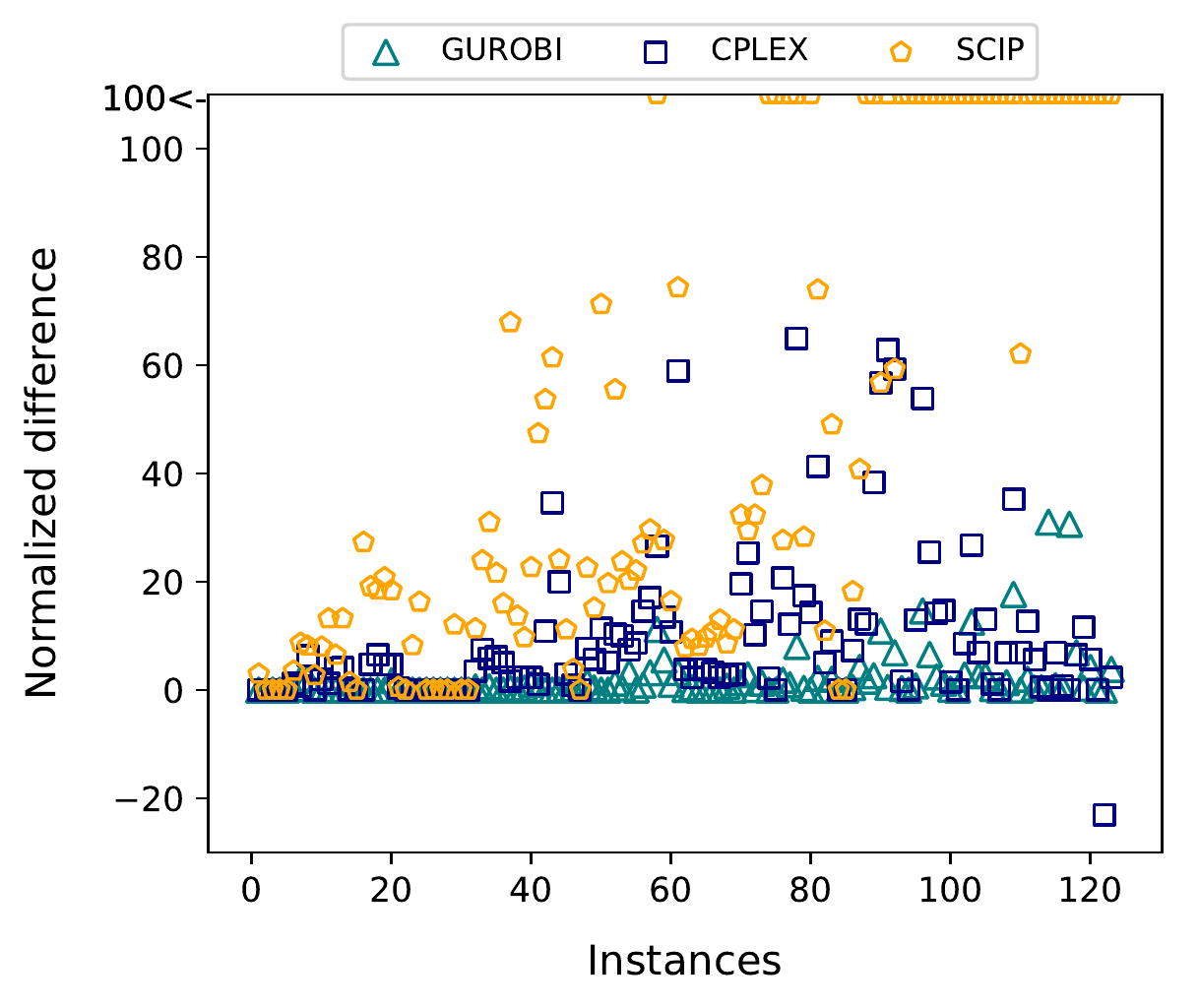}
		\caption{IP solvers}
		\label{subfig:QAP_IPsolvers_120sec}
	\end{subfigure}
	\begin{subfigure}{0.49\textwidth}
		\centering
        \includegraphics[scale=0.62]{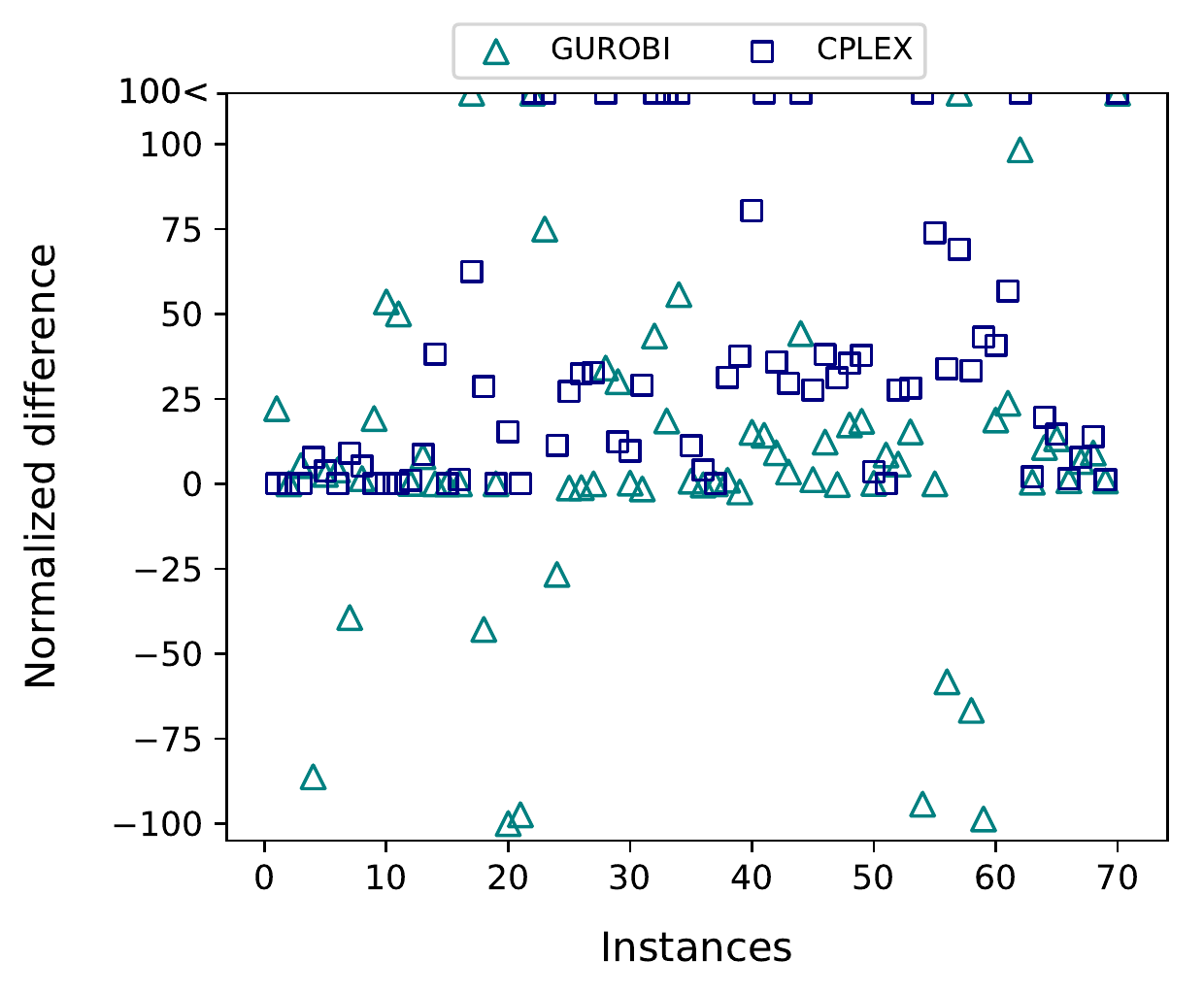}
        \caption{QUBO solvers}
        \label{subfig:QAP_QUBOsolvers_exact_120sec}
	\end{subfigure}
	\caption{Performance comparison of GUROBI, CPLEX and SCIP for QAP 
	under 120-second time limit, with 10-minute results of GUROBI as IP solver taken as reference.}
	\label{fig:QAP_exact_comp}
\end{figure}

{\partitlestyle Comparison of exact solvers.} 
Next, we compare GUROBI, CPLEX and SCIP to each other in Figure \ref{fig:QAP_exact_comp}. 
In this case, we take GUROBI's 10-minute performance (as IP solver) as reference and report the normalized differences with respect to the objective values thereof.
Figure \ref{subfig:QAP_IPsolvers_120sec} shows the results from the three solvers in solving the IP formulations under 120-second time limit. 
We see that GUROBI performs better than CPLEX and SCIP in almost all the instances, and SCIP performs poorly especially for large instances, as revealed by the data points with a normalized difference value of larger than 100. 
In comparing the performances in solving the QUBO formulations, we exclude SCIP, because it can barely find any feasible solutions (as reported in Table \ref{tab:QAP_feasibility}).
The results from GUROBI and CPLEX under 120-second time limit provided in Figure \ref{subfig:QAP_QUBOsolvers_exact_120sec} indicate that GUROBI mostly outperforms CPLEX in solving the QUBO formulations of the QAP instances. 
Additionally, for a subset of instances, GUROBI as QUBO solver provides better solutions than it does as  IP solver.


\begin{figure}[ht]
    \centering
	\begin{subfigure}{0.49\textwidth}
		\centering
	    \includegraphics[scale=0.62]{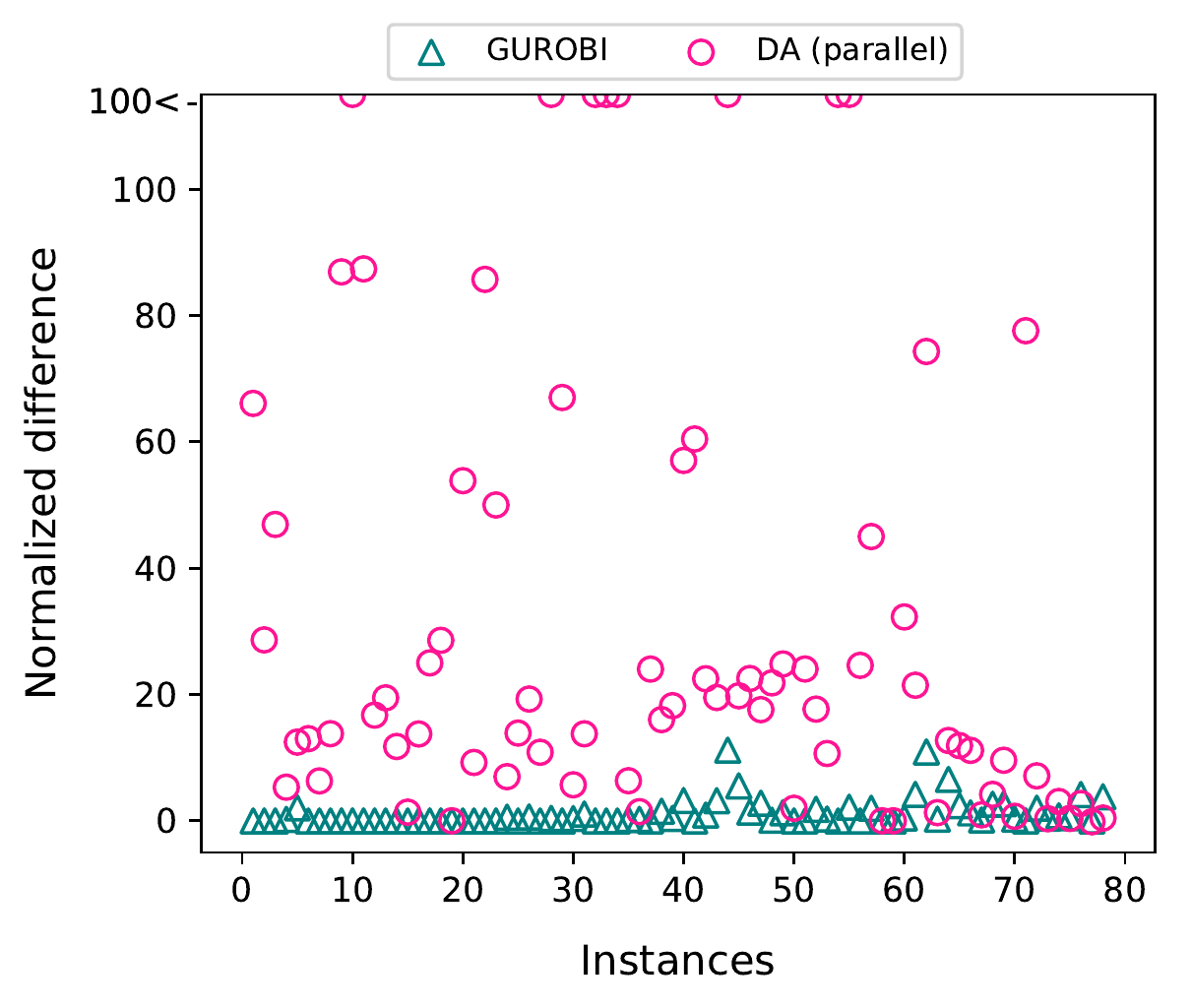}
		\caption{GUROBI as IP solver}
		\label{subfig:QAP_gap_da_vs_gurobi_IP}
	\end{subfigure}
	\begin{subfigure}{0.49\textwidth}
		\centering
        \includegraphics[scale=0.62]{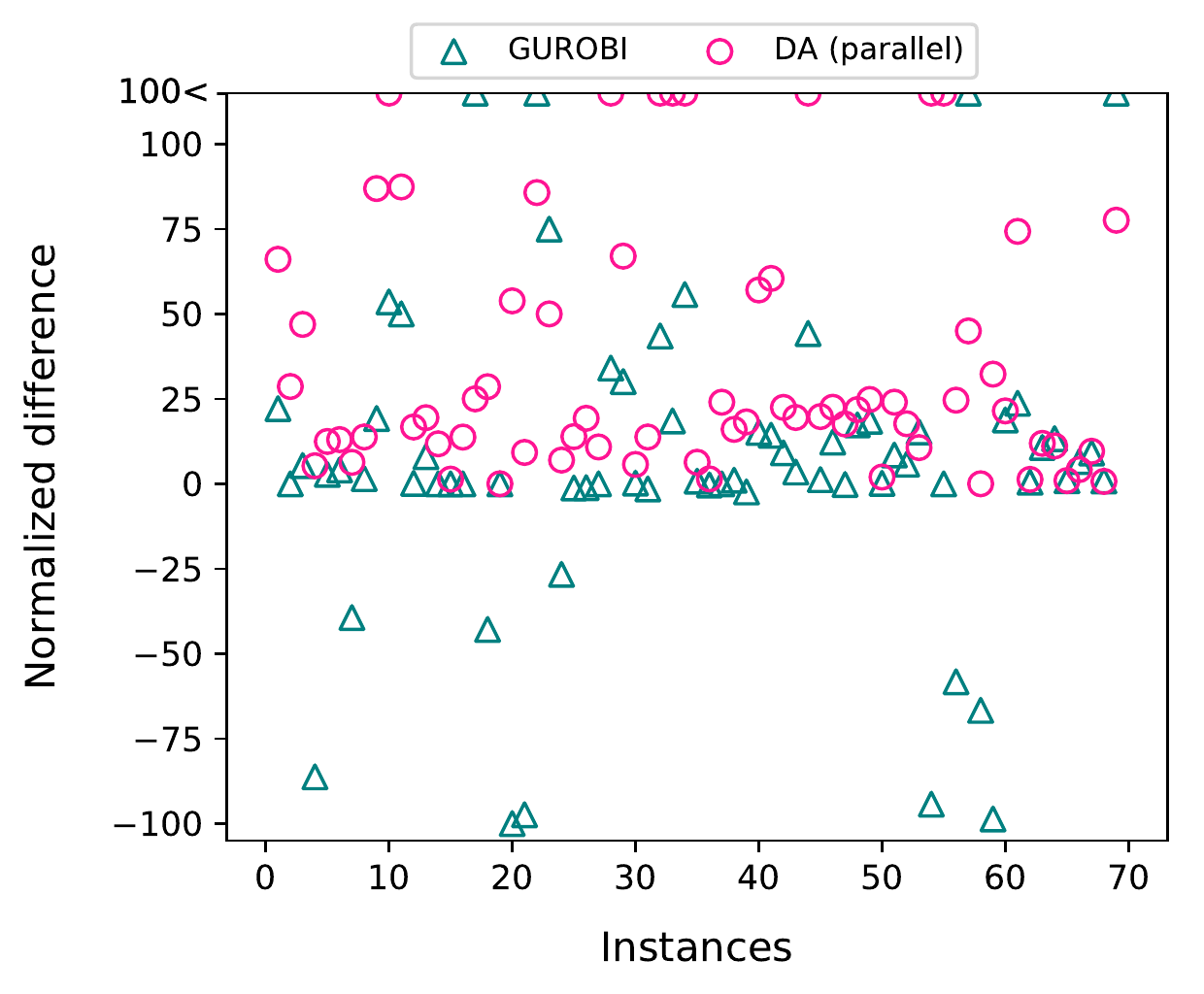}
        \caption{GUROBI as QUBO solver}
        \label{subfig:QAP_gap_da_vs_gurobi_QUBO}
	\end{subfigure}
	\caption{Performance comparison of DA and GUROBI for QAP under 120-second time limit, with 10 minute results of GUROBI as IP solver taken as reference.}
	\label{fig:QAP_da_vs gurobi}
\end{figure}

{\partitlestyle Comparison of DA and exact solvers.} 
Based on the above results, we choose GUROBI and parallel mode of DA for further evaluation, and compare them in Figure \ref{fig:QAP_da_vs gurobi} using their 120-second performances.
From Figures \ref{subfig:QAP_gap_da_vs_gurobi_IP} and \ref{subfig:QAP_gap_da_vs_gurobi_QUBO}, we see that GUROBI as IP solver always yields at least as good or better objective values than DA, where DA occasionally does better than GUROBI as QUBO solver.
As we take GUROBI's 10-minute results as reference, Figure \ref{subfig:QAP_gap_da_vs_gurobi_IP} also implies that GUROBI cannot really improve the quality of the solutions with increased time limit, because the normalized differences are almost always zero. 



{\partitlestyle Overall performances.} 
We conclude our experimental analysis for QAP with an overall summary of the computational results, provided in Table \ref{tab:qap_summary}. 
The table reports the average normalized differences with 10 minute results of GUROBI as IP solver taken as reference (``Norm diff''), average run times (``Time (sec)''), and average percentage gap provided by each solver (``\% solver gap''), for DA and the other three as IP solvers. 
All of the reported results are calculated for a set of common instances for which all solvers can provide feasible solutions in the given time limits. 
As before, the \largeVal~signs in the table stand for values that are greater than $10^5$.
The average performance of GUROBI is better than the other exact solvers and DA in terms of solution quality. 
In terms of run times, on the other hand, DA in normal mode is able to provide feasible solutions in less than half a second. 

\input{tables/QAP_summary}

%% file: tables/QAP_feasibilityTable.tex
\begin{table}[ht]
	\centering
	\caption{Feasibility analysis of DA and exact solvers under 120-second time limit for QAP.}
	\label{tab:QAP_feasibility}
	 \def\arraystretch{1}
    \scalebox{0.84}{
	\begin{tabular}{crrrrrrrr}
		\toprule
		& \multicolumn{3}{c}{IP} & \multicolumn{5}{c}{QUBO} \\
		\cmidrule(lr){2-4}\cmidrule(lr){5-9}
		\parbox[t]{0.5cm}{\centering \phantom{o}} 
		& \parbox[t]{1.6cm}{\centering GUROBI} 
		& \parbox[t]{1.3cm}{\centering CPLEX} 
		& \parbox[t]{1.3cm}{\centering SCIP}
		& \parbox[t]{1.3cm}{\centering DA (n)} 
		& \parbox[t]{1.3cm}{\centering DA (p)}
		& \parbox[t]{1.6cm}{\centering GUROBI} 
		& \parbox[t]{1.3cm}{\centering CPLEX} 
		& \parbox[t]{1.3cm}{\centering SCIP}  \\
		\midrule
		\% feas & 100.0 & 100.0 & 100.0 & 63.4 & 63.4  & 61.0 & 57.7 & 6.5  \\[0.075cm]
		\% solver gap  & 84.7 & 88.5 & 93.8 & {--} & {--} & {--} & {--} & {--}   \\[-0.05cm]
		(DA infeas)\\
		\bottomrule
	\end{tabular}%
	}
\end{table}%

%% file: tables/QAP_summary.tex
\begin{table}[htbp]
  \centering
  \caption{Average performance of DA, GUROBI, CPLEX and SCIP as IP solvers for QAP.}
  \def\arraystretch{1}
  \scalebox{0.75}{
    \begin{tabular}{crrrrrrrrrrr}
    \toprule
      &  & \multicolumn{3}{c}{60-sec limit} & \multicolumn{4}{c}{120-sec limit} & \multicolumn{3}{c}{10-min limit} \\
      \cmidrule(lr){3-5}\cmidrule(lr){6-9}\cmidrule(lr){10-12}
      & \parbox[t]{1.3cm}{\centering DA (n)} 
      & \parbox[t]{1.4cm}{\centering GUROBI} 
      & \parbox[t]{1.1cm}{\centering CPLEX}  
      & \parbox[t]{1.1cm}{\centering SCIP} 
      & \parbox[t]{1.3cm}{\centering DA (p)}  
      & \parbox[t]{1.4cm}{\centering GUROBI} 
      & \parbox[t]{1.1cm}{\centering CPLEX} 
      & \parbox[t]{1.1cm}{\centering SCIP} 
      & \parbox[t]{1.4cm}{\centering GUROBI} 
      & \parbox[t]{1.1cm}{\centering CPLEX} 
      & \parbox[t]{1.1cm}{\centering SCIP} \\
    \midrule
    Norm diff & 66.7 & 1.9 & 9.9 & \largeVal & 37.5 & 1.0 & 8.0 &  \largeVal & 0.0 & 5.2 & \largeVal \\[0.05cm]
    Time (sec) & 0.4  & 46.4 & 49.3 & 57.3 & 106.8 & 90.3 & 95.5 & 112.1 & 417.1 & 460.5 & 509.3\\[0.05cm]
    \% solver gap & {--} & 63.7 & 71.3 & 91.1 & {--} & 59.4 & 69.6 & 86.8 & 53.5 & 63.6 & 76.9\\[0.05cm]
    \bottomrule
    \end{tabular}%
    }
  \label{tab:qap_summary}%
\end{table}%


%% file: QCPP.tex
\subsection{Quadratic cycle partition problem}
\label{subsec:QCPP}

In this section, we explore a graph application, which, to the best of our knowledge, has not been addressed with DA or similar solution technologies before.
Given a graph $G=(V,A)$ with vertex set $V$ and arc set $A$, a {\it cycle partition} 
is a set of (directed) cycles such that every vertex is contained in exactly one cycle.
The {\it cycle partition problem} aims to find a cycle partition of minimum (or maximum) weight, and in that regard the well-known travelling salesperson problem is a special case of it.
The cycle partition problem can be defined on undirected graphs as well, but in this study, we consider the problem on directed graphs.
In particular, we consider the {\it quadratic cycle partition problem} (QCPP), which aims to find a cycle partition in a directed graph so that the sum of the interaction costs between consecutive arcs is minimized.

QCPP is an important problem in that it is closely related to the quadratic traveling salesperson problem \citep{jager2008algorithms}, 
and has many variants with applications in various fields \citep{meijer2020quadratic}.
Different solution approaches for QCPP or its variants have been studied in the literature, such as local search methods, column generation to compute lower bounds, and branch-and-bound algorithms  \citep{jager2008algorithms,galbiati2014minimum}.
In the sequel, we test the IP, QUBO, and DA approaches on QCPP.

\subsubsection{Formulations}

QCPP can be formulated as an IP model containing only binary variables and equality constraints. 
Suppose that we are given a directed graph $G=(V,A)$ with vertex set $V = \{1,...,n\}$ and arc set $A$, and arc pair costs $q_{{a_1}{a_2}}$'s for all $a_1, a_2 \in A$.
The cost $q_{{a_1}{a_2}}$ of an arc pair $\{a_1, a_2\}$ is nonzero only if the associated arcs are consecutive, that is, if the vertex of arc $a_1 \in A$ is directed to and the vertex of arc $a_2 \in A$ is directed from are the same.
Let $\delta^{+}(i)$ and $\delta^{-}(i)$ respectively denote the set of all outgoing and incoming arcs of vertex $i \in V$. 
We define binary decision variables $x_{a}$'s 
for the assignment of arcs to cycles as such:
\begin{equation*}
    x_a = 
    \begin{cases}
        1, & \text{if arc $a$ is assigned to a cycle} \\
      0, & \text{otherwise}.
    \end{cases}
\end{equation*}
Then, an IP formulation for QCPP can be written as follows:
\begin{subequations}
\begin{alignat}{2}
	\min\quad &{\sum_{a_1 \in A}\sum_{a_2 \in A} 
	q_{{a_1}{a_2}} x_{a_1} x_{a_2}} \label{eq:QCPP_IP_obj}\\
	\text{s.t.} \quad &{ \sum_{a \in \delta^{+}(i)} x_{a} = 1} \qquad && { i \in V } \label{eq:QCPP_IP_c1} \\[0.1cm]
	&{ \sum_{a \in \delta^{-}(i)} x_{a} = 1} && {i \in V} \label{eq:QCPP_IP_c2}\\[0.1cm]
	& x_{a} \in \{0,1\} && {a \in A} \label{eq:QCPP_IP_c3}
\end{alignat}
\label{m:QCPP}
\end{subequations}%
\indent
The quadratic objective function \eqref{eq:QCPP_IP_obj} minimizes the total interaction cost between selected arcs. Constraints \eqref{eq:QCPP_IP_c1} and \eqref{eq:QCPP_IP_c2} ensure that each vertex must be incident to exactly one cycle. 
Constraints \eqref{eq:QCPP_IP_c3} restrict the decision variables to be binary.

Similar to QAP, we can dualize the equality constraints by taking the square of the difference between the left- and right-hand sides of each constraint and adding them to the objective as penalty terms.
Then, a QUBO formulation for QCPP can be written as follows:
\begin{subequations}
\label{m:QCPP_qubo}
\begin{align}
\min\quad &{ \sum_{a_1 \in A} \sum_{a_2 \in A} q_{{a_1}{a_2}} x_{a_1} x_{a_2} } + \penaltyCoef \ \sum_{i \in V}\left(\sum_{a \in \delta^{+}(i)} x_{a} - 1 \right)^2 +  \penaltyCoef  \sum_{i \in V}\left(\sum_{ a \in \delta^{-}(i)} x_{a} - 1 \right)^2  \label{eq:QCPP_qubo_obj}\\[0.25cm]
\text{s.t.} \quad & x_{a} \in \{0,1\} \qquad  a \in A, \label{eq:QCPP_qubo_c}
\end{align}
\end{subequations}
\noindent where $\penaltyCoef > 0$ is the penalty coefficient. 
For a sufficiently large value of $\penaltyCoef$, the QUBO formulation \eqref{m:QCPP_qubo} produces provably optimal solutions for the original IP model \eqref{m:QCPP}.

\subsubsection{Instances}
We generated 80 QCPP instances of varying sizes by following the procedure in \citep{meijer2020sdp}.
Specifically, we first created the graph instances according to the \cite{erdos1959random} model, where each possible arc is added with a certain probability.
For each graph instance, the interaction cost between any pair of successive arcs is chosen as an integer uniformly distributed between 0 and 100.
The number of vertices in our instances vary between 25 to 175, and the arc densities range from 0.25 to 0.75.
We note that the instances used in \citep{meijer2020sdp} are relatively small, and typically easy to handle by the state-of-the-art solvers now.
Thus, we created our test bed by including larger and hence more challenging QCPP instances.
We provide brief information on our test bed in Table \ref{tab:qcc_inst_info}.

\input{tables/QCPP_instance_info}

In order to ensure that our model correctly represents the problem for each instance, we made sure that every vertex of a graph instance has at least one incoming and one outgoing arc, because otherwise there will be at least one vertex for which constraints of type \eqref{eq:QCPP_IP_c1} and \eqref{eq:QCPP_IP_c2} are incomplete or absent in the model, which would lead to an incorrect representation of the problem.
Even though we ensure that all the required constraints are contained the IP model, it is not trivial to guarantee the existence of a cycle partition in a graph, mainly because every vertex needs to be covered by exactly one cycle.
Nonetheless, all the problem instances in our test bed have a nonempty feasible solution space.

\subsubsection{Experimental results}
In this section, we present the results of our computational study for QCPP,
starting with pictorial expositions and concluding with tabular summaries.
The points on the $x$-axes of figures represent individual instances arranged in ascending order with respect to the number of variables they have.
We use normalized differences as a means to compare solvers, where the objective values (upper bounds) from GUROBI's 10-minute experiments are taken as reference. 

\begin{figure}[ht]
    \centering
	\begin{subfigure}{0.49\textwidth}
		\centering
	    \includegraphics[scale=0.62]{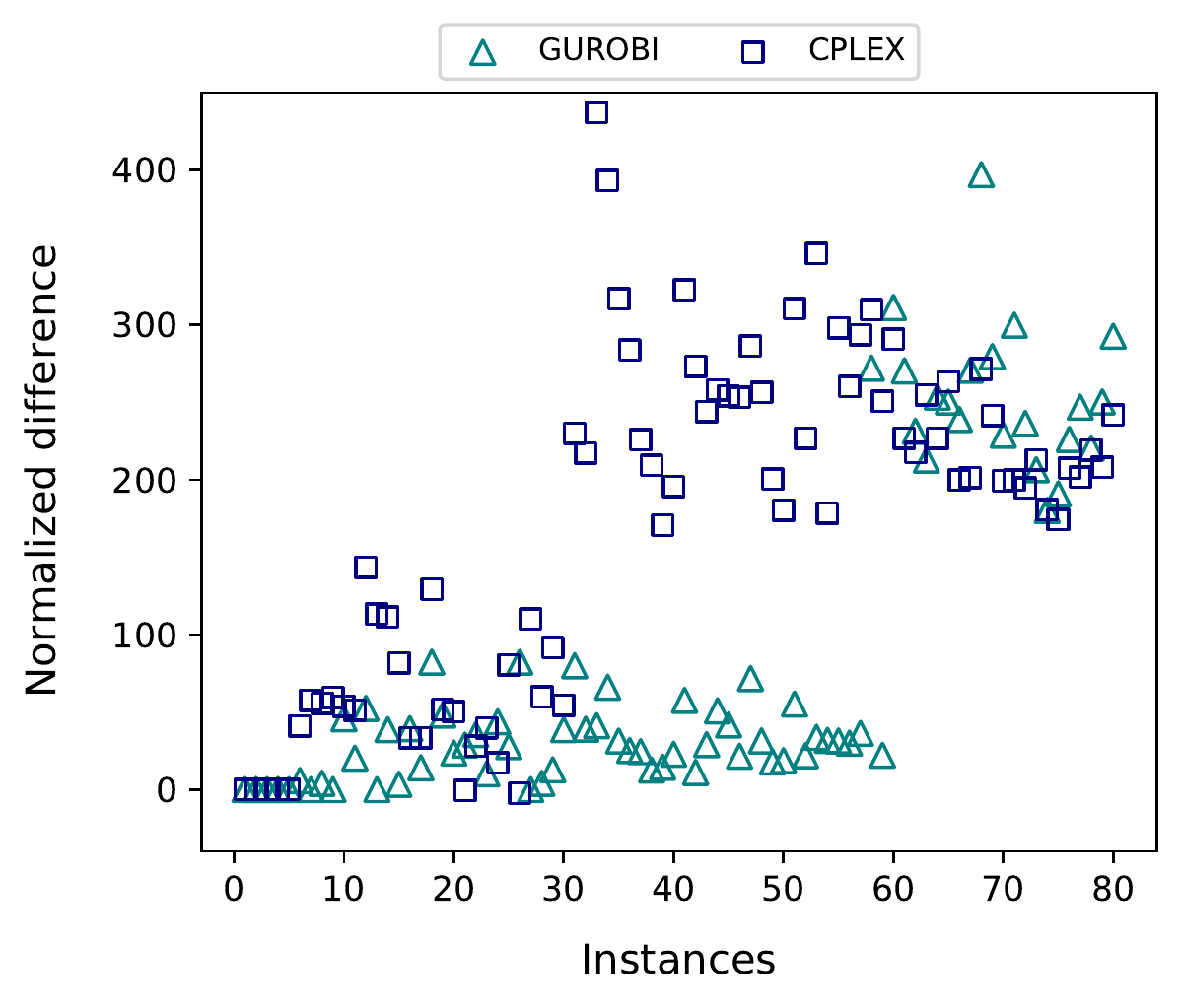}
		\caption{60-second time limit}
		\label{subfig:qcc_cplex_vs_gurobi_normUBdiff_60}
	\end{subfigure}
	\begin{subfigure}{0.49\textwidth}
		\centering
	    \includegraphics[scale=0.62]{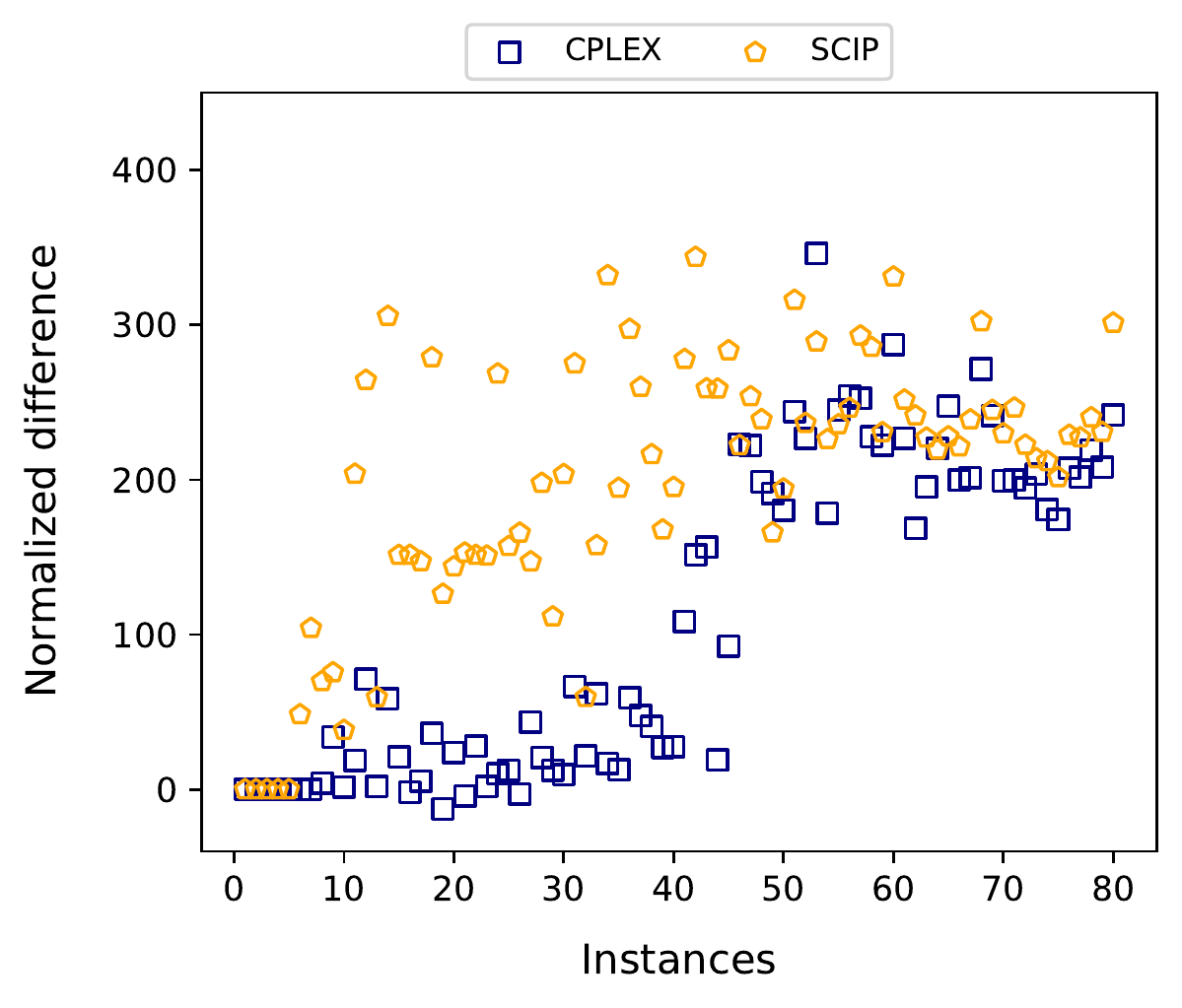}
		\caption{10-minute time limit}
		\label{subfig:qcc_cplex_vs_scip_normUBdiff_600}
	\end{subfigure}
	\caption{Performance comparison of GUROBI, CPLEX and SCIP for QCPP, with 10-minute results of GUROBI as IP solver taken as reference.}
	\label{fig:qcc_cplex_vs_gurobi_vs_scip_normUBdiff}
\end{figure}

{\partitlestyle Comparison of exact solvers.} 
We first compare the performances of GUROBI, CPLEX and SCIP in solving the IP formulations in terms of normalized differences, shown in Figure \ref{fig:qcc_cplex_vs_gurobi_vs_scip_normUBdiff}.

Figure \ref{subfig:qcc_cplex_vs_gurobi_normUBdiff_60} contains the normalized differences of 60-second results from GUROBI and CPLEX. 
In this case, we do not include SCIP, because it failed to deliver any valid lower or upper bounds, even when the time limit is increased to 120 seconds. Our observation is that GUROBI does better than CPLEX especially in small to moderate-size instances.
Specifically, GUROBI does at least as good as CPLEX in 72.5\% of the instances.

We also compare the performance of solvers under 10-minute time limit, shown in Figure \ref{subfig:qcc_cplex_vs_scip_normUBdiff_600}.
Here, we do not explicitly show the results from GUROBI, because the associated upper bound values are taken as reference and hence the corresponding normalized difference values would simply be anchored at zero level.
We see from this figure that both CPLEX and SCIP almost always yield inferior upper bounds as compared to GUROBI, because the data points are almost always above zero.
While CPLEX and SCIP perform the same on small instances, CPLEX does typically better than SCIP in the rest, and the difference is especially noticable in moderate-size instances.
In particular, CPLEX does at least as well as SCIP in 91.3\% of the instances.
As a result of this comparative evaluation, we continue with GUROBI to compare DA to, as it yields the most promising results among the three solvers.

\begin{figure}[ht]
    \centering
	\begin{subfigure}{0.49\textwidth}
		\centering
	    \includegraphics[scale=0.62]{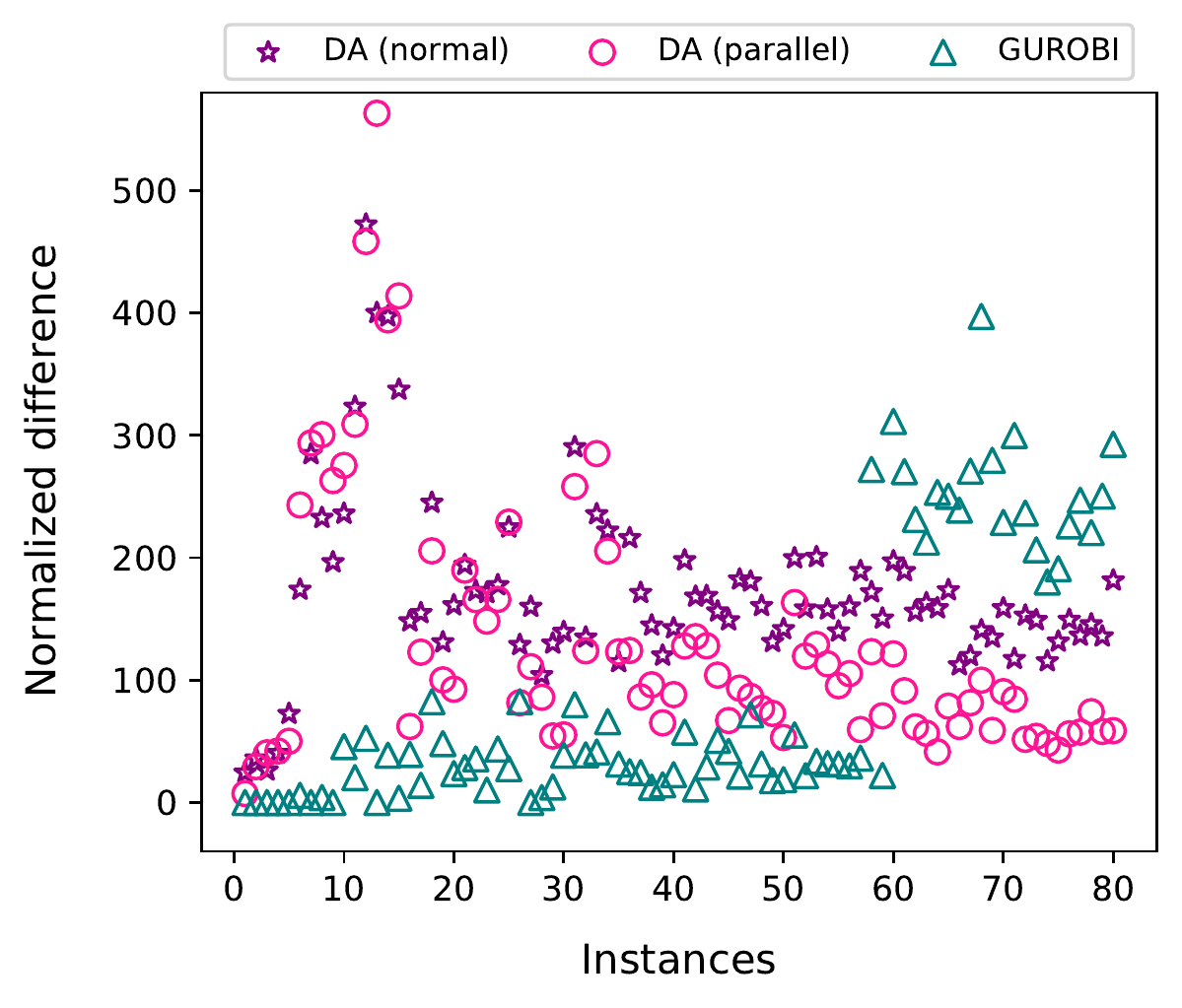}
		\caption{60-second time limit}
		\label{subfig:qcc_da_vs_gurobi_normUBdiff_60}
	\end{subfigure}
	\begin{subfigure}{0.49\textwidth}
		\centering
	    \includegraphics[scale=0.62]{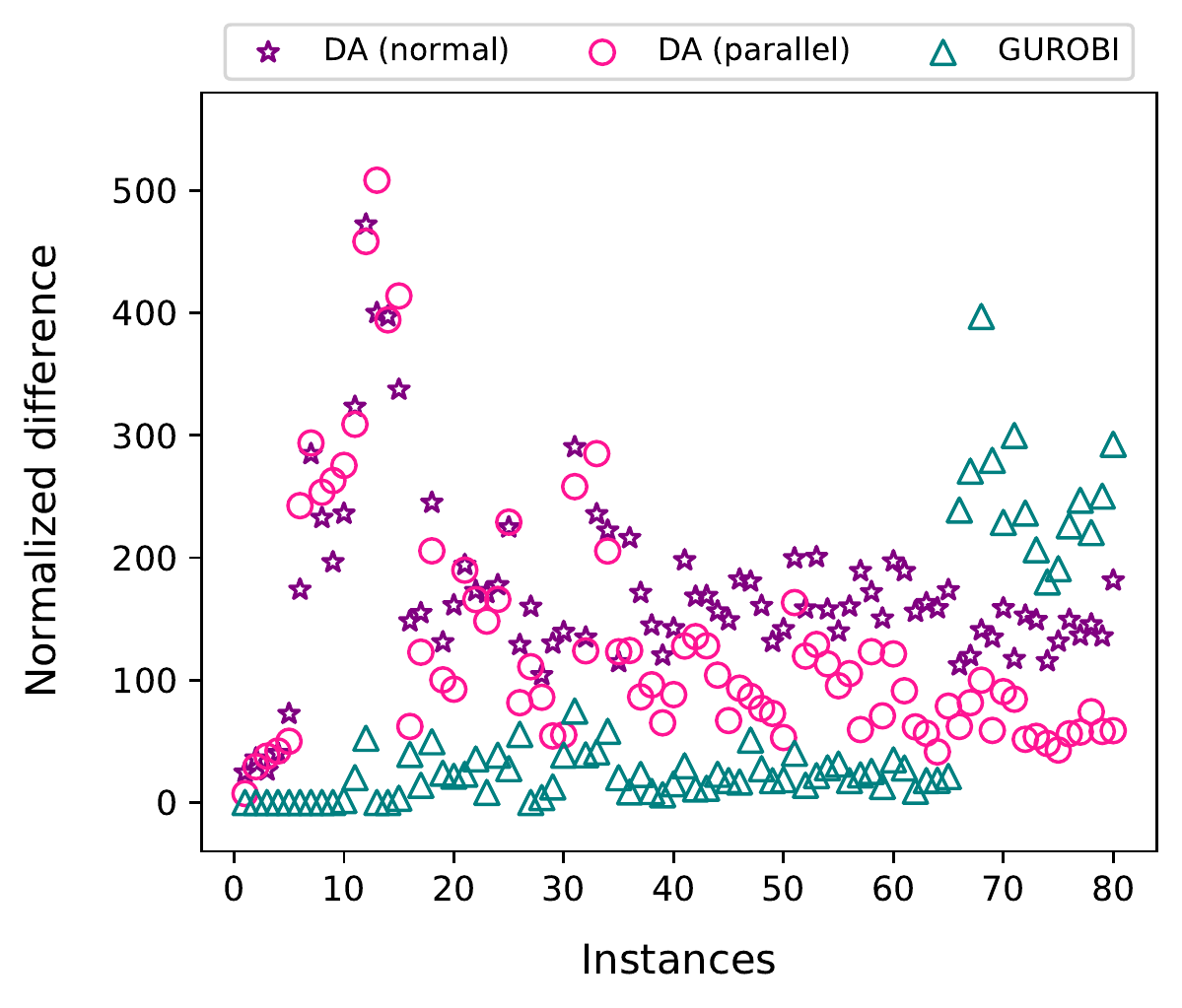}
		\caption{120-second time limit}
		\label{subfig:qcc_da_vs_gurobi_normUBdiff_120}
	\end{subfigure}
	\caption{Performance comparison of DA and GUROBI for QCPP, with 10-minute results of GUROBI as IP solver taken as reference.}
	\label{fig:qcc_da_vs_gurobi_normUBdiff}
\end{figure}

{\partitlestyle Comparison of DA and exact solvers.} 
Figure \ref{fig:qcc_da_vs_gurobi_normUBdiff} shows the normalized differences for GUROBI and both modes of DA.
We see that GUROBI yields better objective values than DA for relatively small instances.
On the other hand, for larger instances, i.e., when we look at the data points towards the right end of the $x$-axis, both normal and parallel modes of DA yield significantly better objective values, with the parallel mode being the best among all.
Particularly, for 60- and 120-second results, DA outperforms in approximately 29\% and 19\% of the instances with at least 4500 and 7100 variables, respectively.
So, the gap between the data points corresponding to GUROBI and the parallel mode of DA is greater for larger instances, and that DA shows a robust performance despite increasing instance sizes.
Therefore, we conclude that DA is particularly successful in yielding good-quality solutions quickly for large instances.

\begin{figure}[H]
    \centering
    \includegraphics[scale=0.62]{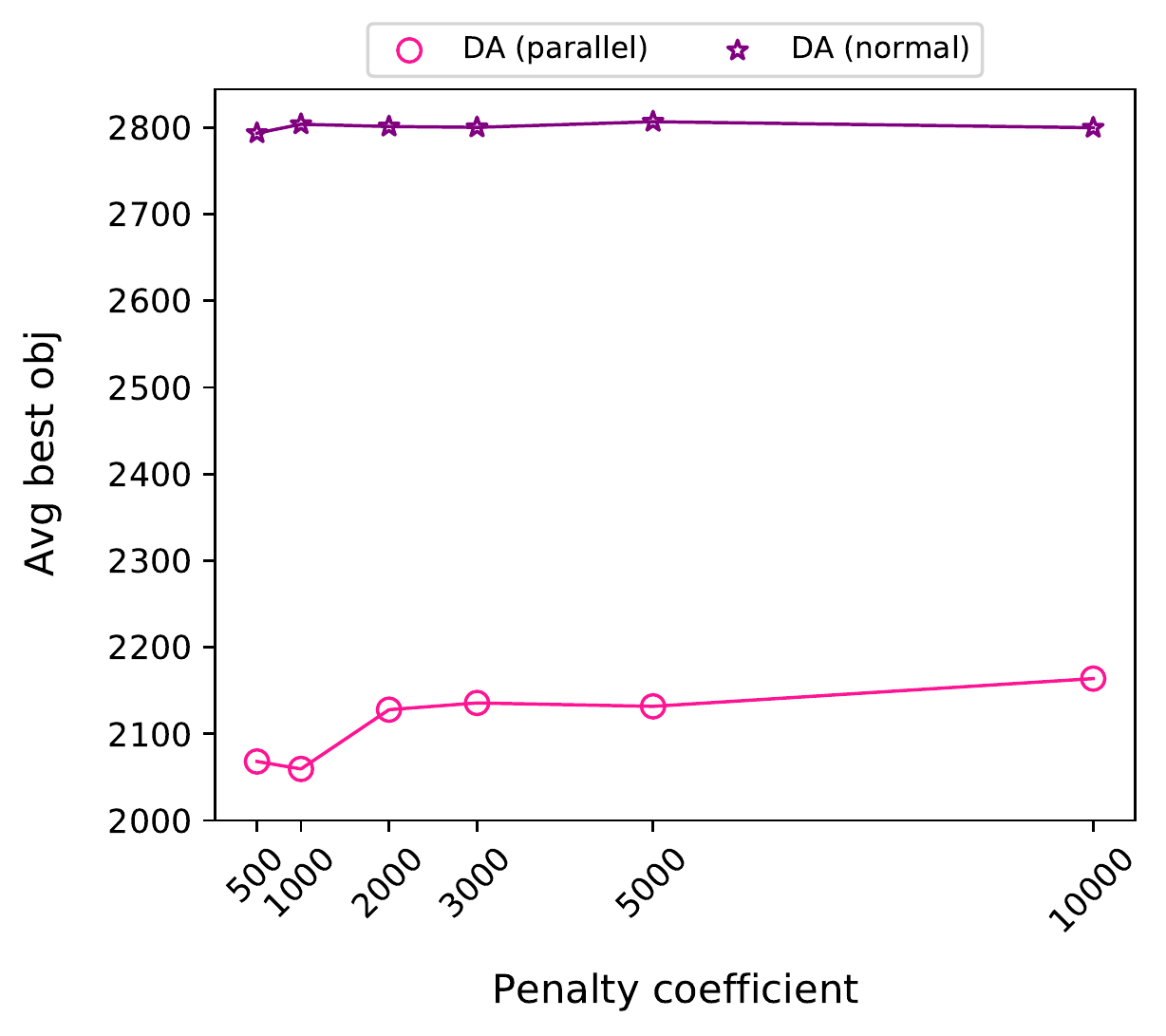}
	\caption{Effect of penalty coefficient values on DA's performance for QCPP.}
	\label{fig:qcc_pentrials}
\end{figure}

{\partitlestyle Penalty coefficient analysis.} 
Next, we investigate the effect of the penalty coefficient values on DA’s performance. 
Figure \ref{fig:qcc_pentrials} shows the average of the best objective values found by DA using six different penalty coefficient values.
The lower end of the values are determined in such a way that DA is able to yield a feasible solution to the entire set of instances we experiment with, and the largest value is set so that the general trend in the performances can be observed.

We first see that the parallel mode of DA delivers significantly better results than the normal mode on the average.
In addition, setting $\penaltyCoef = 1000$ gives the least overall average values for parallel mode, and the performance keeps deteriorating for larger values of the penalty coefficient $\penaltyCoef$.
As the parallel mode yields better results in general, we took its minimum point as a reference and used $\penaltyCoef = 1000$ when running the QUBO models for QCPP, i.e., in obtaining the above-presented results.
We note that even if we were to carry out our comparative analysis with a different penalty coefficient, our conclusions would not generally change, because the difference it would create in DA's performance is less than how much it differs from GUROBI for cases where DA outperforms.

{\partitlestyle Overall performances.} 
In the remainder of this section, we provide tabular summaries of overall performances.
We start with the average performance of solving the QUBO formulations with DA in comparison to solving the IP formulations with GUROBI, CPLEX and SCIP under different time limits, reported in Table \ref{tab:qcc_summary}.
The table contains three blocks of columns reserved for different time limits, and another column dedicated to the normal mode of DA only, because its execution times do not exceed half a second even when the number of iterations is increased greatly. 
The blocks contain the results of solving the QUBO formulations with DA and/or solving the IP formulations with GUROBI, CPLEX and SCIP.
For each option (column), we report the averages of best found objective values (``UB"), normalized differences with 10 minute results of GUROBI as IP solver taken as reference (``Norm diff"), and solution times (``Time (sec)").
We note that the normalized differences are calculated by taking the objective values from 10-minute experiments with GUROBI as reference, because it yielded the best values.
Furthermore, SCIP is excluded from 60-second and 120-second blocks, because it failed to yield any valid lower or upper bound within these time limits.

\input{tables/QCPP_summary_IP}

Under 60-second and 120-second time limits, parallel mode of DA gives the best average upper bounds, and the difference from other solvers is considerably large in the 60-second case.
As for the normal mode, its average upper bound value is better than those obtained from CPLEX and SCIP under 10-minute limit, and its execution time is only 0.4 second. 
In terms of the average normalized differences, GUROBI averages are the best under all three time limits.
Combining these aggregated results with those over the individual instances discussed above, we can say that DA is particularly successful in yielding good-quality solutions for relatively large instances, and if such solutions are needed within a second, normal mode of DA is the most promising option.
We note that all the solutions that DA yield are feasible, and hence is not explicitly reported in the table.

We finally evaluate the overall performances of GUROBI, CPLEX and SCIP in solving the QUBO formulations, and report the associated results in Table \ref{tab:qcc_summary_qubo}.
In this case, since we are solving QUBO models, we include the feasibility percentage of found solutions (``\% feas") in addition to the three measures used in Table \ref{tab:qcc_summary}.
We show values that are greater than $10^5$ with the \largeVal~sign in the table.

\input{tables/QCPP_summary_QUBO}

We see that SCIP fails to deliver any feasible solution within one minute, and in two minutes, it achieves feasibility in only 6.3\% of the instances.
Even though GUROBI and CPLEX do significantly better relative to SCIP, their performance as QUBO solvers is still not satisfactory; the feasibility percentage goes only as much as 67.5\% in the best case. 
So, we conclude that DA is the best QUBO solver, and yields the best results for large-scale instances in cases where short solution times are sought.

%% file: tables/QCPP_instance_info.tex
\begin{table}[htbp]
  \centering
  \caption{QCPP instance information.}
   \def\arraystretch{1}
  \scalebox{0.86}{
    \begin{tabular}{cccc}
    \toprule
      \parbox[t]{2.2cm}{\centering \# instances} & \parbox[t]{2.2cm}{\centering \# vertices} & \parbox[t]{2.8cm}{\centering Density} & \parbox[t]{2.2cm}{\centering \# variables} \\
    \midrule
     80    & 25--175 & 0.25, 0.50, 0.75 & 127--7963 \\    
     \bottomrule
    \end{tabular}%
    }
  \label{tab:qcc_inst_info}%
\end{table}%

%% file: tables/QCPP_summary_IP.tex
\begin{table}[htbp]
  \centering
  \caption{Average performance of DA, GUROBI, CPLEX and SCIP as IP solvers for QCPP.}
  \def\arraystretch{1}
  \scalebox{0.77}{
    \begin{tabular}{crrrrrrrrrr}
    \toprule
      &       & \multicolumn{3}{c}{60-sec limit} & \multicolumn{3}{c}{120-sec limit} & \multicolumn{3}{c}{10-min limit} \\
      \cmidrule(lr){3-5}\cmidrule(lr){6-8}\cmidrule(lr){9-11}
      & \parbox[t]{1.3cm}{\centering DA (n)} 
      & \parbox[t]{1.3cm}{\centering DA (p)} 
      & \parbox[t]{1.5cm}{\centering GUROBI} 
      & \parbox[t]{1.3cm}{\centering CPLEX}  
      & \parbox[t]{1.3cm}{\centering DA (p)}  
      & \parbox[t]{1.5cm}{\centering GUROBI} 
      & \parbox[t]{1.3cm}{\centering CPLEX} 
      & \parbox[t]{1.5cm}{\centering GUROBI} 
      & \parbox[t]{1.3cm}{\centering CPLEX} 
      & \parbox[t]{1.3cm}{\centering SCIP} \\
    \midrule
    UB & 2803.7 & 2059.5 & 2573.5 & 3311.9 & 2057.6 & 2152.6 & 3222.1 & 1105.7 & 2812.9 & 3522.4 \\[0.05cm]
    Norm diff & 172.1 & 129.5 & 90.0 & 173.3 & 128.2 & 63.8 & 161.6 & 0.0  & 112.6 & 199.3 \\[0.05cm]
    Time (sec) & 0.4  & 61.7 & 56.3 & 57.4 & 103.8 & 112.0 & 114.4 & 541.0 & 563.7 & 567.1 \\[0.05cm]
    \% gap & 91.8 & 91.8 & 84.0 & 86.7 & 91.7 & 83.3 & 85.8 & 82.2 & 83.5 & 88.0 \\[0.05cm]
    \% solver gap & {--} & {--} & 87.1 & 93.0 & {--} & 85.0 & 92.5 & 82.2 & 90.5 & 93.5 \\[0.05cm]
    \bottomrule
    \end{tabular}%
    }
  \label{tab:qcc_summary}%
\end{table}%

%% file: tables/QCPP_summary_QUBO.tex
\begin{table}[htbp]
  \centering
  \caption{Average performance of GUROBI, CPLEX and SCIP as QUBO solvers for QCPP.}
  \def\arraystretch{1}
  \scalebox{0.83}{
    \begin{tabular}{crrrrrr}
    \toprule
      & \multicolumn{3}{c}{60-sec limit} & \multicolumn{3}{c}{120-sec limit} \\
      \cmidrule(lr){2-4}\cmidrule(lr){5-7}
      & \parbox[t]{2cm}{\centering GUROBI} & \parbox[t]{2cm}{\centering CPLEX} & \parbox[t]{2cm}{\centering SCIP} & \parbox[t]{2cm}{\centering GUROBI} & \parbox[t]{2cm}{\centering CPLEX} & \parbox[t]{2cm}{\centering SCIP} \\
    \midrule
    UB & 91802.8 & 6143.7 & 158424.1 & 77820.2 & 5967.7 & 157672.9 \\[0.05cm]
    Norm diff & 5701.4 & 453.3 & \largeVal & 4106.9 & 416.5 & \largeVal \\[0.05cm]
    Time (sec) & 60.8 & 103.3 & 60.4 & 120.2 & 138.6 & 121.1 \\[0.05cm]
    \% feas & 55.0  & 30.0  & 0.0  & 67.5  & 32.5  & 6.3 \\[0.05cm]
    \% gap & 86.5 & 89.3 & 99.5 & 86.2 & 87.5 & 96.4 \\[0.05cm]
    \% solver gap & 2507.5 & \largeVal & \largeVal & 2068.5 & \largeVal & \largeVal  \\[0.05cm]
    \bottomrule
    \end{tabular}%
    }
  \label{tab:qcc_summary_qubo}%
\end{table}%

%% file: SelCol.tex
\subsection{Selective graph coloring problem} \label{subsec:selcol}

In this section, we focus on the selective graph coloring problem, Sel-Col, which is a generalization of the classical graph coloring problem. 
Coloring of a graph is assignment of labels (colors) to its vertices, such that no pair of vertices that are linked by an edge receives the same color.
Given a graph and a partition of its vertex set into {\it clusters}, the aim in Sel-Col is to choose exactly one vertex per cluster so that, among all possible selections, the number of colors necessary to color the selected set of vertices is minimized. 
As opposed to the classical graph coloring problem, we do not color every vertex in this problem; we only do so for the selected set.
Figure \ref{fig:selcol_example} illustrates a small Sel-Col instance together with one feasible and one optimal selective coloring of it. 
To the best of our knowledge, this problem has not been addressed with DA or similar solution technologies before.

\input{figures/SelCol/illustrativeExample}

Sel-Col has emerged as a model to select a route and assign a proper wavelength to each in the second step of a two-phase solution strategy for the routing and wavelength assignment problem that arises in optical networks \citep{li2000partition}.
The problem has many applications in various domains, such as timetabling, constraint encoding, antenna positioning and frequency assignment \citep{demange2015some}.
Several solution approaches have been proposed for Sel-Col, for general instances or specific graph classes, including heuristics \citep{li2000partition}, and exact algorithms such as a branch-and-cut algorithm \citep{frota2010branch}, branch-and-price algorithms \citep{hoshino2011branch,furini2017exact}, and cutting plane algorithms \citep{cseker2019decomposition,cseker2020exact}.
In what follows, we present the results of IP, QUBO and DA approaches in handling Sel-Col.

\subsubsection{Formulations}

Suppose that we are given an undirected graph $G=(V,E)$ with $V = \{1, \ldots ,n\}$, as well as a partition $\mathcal{V}$ of its vertex set into $\mathit{P}$ clusters  $V_1,...,V_P$. 
Sel-Col can be formulated as an IP model, which only contains binary variables but has both equality and inequality constraints. 
We introduce two sets of binary variables. First, we define $y_{k}$'s for $k \in \{1,\ldots,P\}$ to keep track of the number of colors used.
\begin{equation*}
    y_k = 
    \begin{cases}
        1, & \text{if color $k$ is used by some vertex} \\
      0, & \text{otherwise}
    \end{cases}
\end{equation*}
Second, we define $x_{ik}$'s for $i \in V, \ k \in \{1,\ldots,P\}$ in order to relate the color usage to the assignment of colors to vertices. 
\begin{equation*}
    x_{ik} = 
    \begin{cases}
      1, & \text{if vertex $i$ is selected and colored with color $k$} \\
      0, & \text{otherwise}
    \end{cases}
\end{equation*}
Then, an IP formulation for Sel-Col, which appears in \citep{cseker2019decomposition,cseker2020exact}, can be written as follows:
\begin{subequations}
\label{m:selcol}
\begin{alignat}{2}
	\min\quad &{\sum_{k=1}^{\mathit{P}}  y_{k}} \label{eq:selcol_IP_obj}\\
	\text{s.t.} \quad &{ \sum_{i \in V_p} \sum_{k=1}^\mathit{P} x_{ik} = 1} \qquad && { p \in \{1,...,\mathit{P}\} } \label{eq:selcol_IP_c1} \\[0.1cm]
	&{ x_{ik} + x_{jk} \leq 1} && {\{i,j\} \in E , \  k \in  \{1,...,\mathit{P}\}} \label{eq:selcol_IP_c2}\\[0.1cm]
	&{ x_{ik} \leq y_k} && {i \in V, \ k \in  \{1,...,\mathit{P}\}} \label{eq:selcol_IP_c3}\\[0.1cm]
	&y_k \in \{0,1\} && k \in  \{1,...,\mathit{P}\} \label{eq:selcol_IP_c4}\\	%
	& x_{ik} \in \{0,1\} && {i \in V, \ k \in  \{1,...,\mathit{P}\}} \label{eq:selcol_IP_c5}
\end{alignat}
\end{subequations}

The objective function \eqref{eq:selcol_IP_obj} minimizes the number of colors used. 
Constraints \eqref{eq:selcol_IP_c1} guarantee that exactly one vertex is selected per cluster. 
The requirement that adjacent vertices cannot be assigned the same color is met through constraints \eqref{eq:selcol_IP_c2}. 
Constraints \eqref{eq:selcol_IP_c3} link $x$ and $y$ variables by forcing $y_k$ to be one if color $k$ is used by some vertex. 
Finally, \eqref{eq:selcol_IP_c4} and \eqref{eq:selcol_IP_c5} respectively set $y$ and $x$ to be binary variables. 
Note that as we are minimizing the number of colors used, when $x_{ik} = 0 \ \forall i \in V$ for some $k \in  \{1,...,\mathit{P}\}$, then $y_k = 0$ in any optimal solution.   

Despite the presence of inequality constraints, we can reformulate the IP model in \eqref{m:selcol} as a QUBO formulation without any need to introduce additional variables. 
We note that the inequality constraints \eqref{eq:selcol_IP_c2} and \eqref{eq:selcol_IP_c3} can be equivalently written as a pair of quadratic equality constraint sets as such:
\begin{taggedsubequations}{} 
\begin{align}
& x_{ik} x_{jk} = 0 && \{i,j\} \in E, \ k \in  \{1,...,\mathit{P}\}  &  \tag{\ref{eq:selcol_IP_c2} - q} \label{eq:selcol_IP_c2_q} \\[0.1cm]
&(x_{ik} - y_k) x_{ik}  = 0 &&  i \in V, \ k \in  \{1,...,\mathit{P}\}   \tag{\ref{eq:selcol_IP_c3} - q}\label{eq:selcol_IP_c3_q}
\end{align}
\end{taggedsubequations}
%
%
\indent The set of feasible $x$ and $y$ values induced by \eqref{eq:selcol_IP_c2} and \eqref{eq:selcol_IP_c3} are the same as those enforced by \eqref{eq:selcol_IP_c2_q} and \eqref{eq:selcol_IP_c3_q}. Therefore, we can safely replace the inequalities with the equivalent quadratic set of inequalities, and transform the IP model into the following QUBO model:
\begin{subequations}
\label{m:selsol_qubo}
\begin{align}
	\min \    &  {\sum_{k=1}^\mathit{P} \ y_k} 
	+  \penaltyCoef \sum_{p=1}^{P}\left(1-\sum_{i \in V_p} \sum_{k=1}^\mathit{P} x_{ik}\right)^2
	+ \penaltyCoef \sum_{\{i,j\} \in E} \sum_{k=1}^\mathit{P} x_{ik}x_{jk}
	+ \penaltyCoef \sum_{i \in V} \sum_{k=1}^\mathit{P} (x_{ik}-y_k)x_{ik}  \label{eq:selcol_qubo_obj} \\[0.25cm]
	\text{s.t.} \ & y_k \in \{0,1\} \qquad \ k \in  \{1,...,\mathit{P}\} \label{eq:selcol_qubo_c1}\\		
	& x_{ik} \in \{0,1\}  \qquad {i \in V, \ k \in  \{1,...,\mathit{P}\}},	\label{eq:selcol_qubo_c2}
\end{align}
\end{subequations}		

\noindent where $\penaltyCoef$ is the penalty coefficient. 
For $\penaltyCoef > P$, the QUBO formulation for Sel-Col given in \eqref{eq:selcol_qubo_obj}--\eqref{eq:selcol_qubo_c2} is exact; that is, it yields provably optimal solutions for the original problem.

\subsubsection{Size reduction}

For our computational experiments, we consider a large suite of 1200 Sel-Col instances. 
One half of the instances contain perfect graphs, and the other half \cite{erdos1959random} graphs. 
We refer to these two groups as PG and ER instances, respectively.
When we first transformed the IP formulations of these Sel-Col instances into QUBO models, around 37\% of them were size-wise eligible for DA. 
This inclined us to explore ways to reduce the sizes of these instances in order to make a larger proportion fit into DA. 
For this purpose, we design a two-phase heuristic to decrease the number of variables.

Our heuristic is inspired from the fact that Sel-Col naturally consists of two parts; the vertex selection, and the coloring of it. 
In the first phase of our heuristic, we choose one vertex per cluster that has the least number of connections to other clusters. 
By doing so, we aim to decrease the number of edges between the selected set of vertices, and hence potentially the number of colors required (although decreasing the number edges in a graph does not necessarily reduce the number of colors needed). 
In the second phase, we use a heuristic coloring algorithm to color the selected set of vertices, where in each step we color a vertex that has the largest number of already colored neighbors, similar to the DSatur algorithm \citep{brelaz1979new}.
This yields an upper bound on the minimum number of colors needed to color the given vertex selection from the first phase, which also serves as an upper bound on the optimal number of colors of that Sel-Col instance. 
This way, we eliminate the need of defining decision variables with a color index larger than the upper bound found.
This technique can be viewed as a variant of the size reduction ideas proposed by \cite{lewis2017quadratic}, where variables are identified and eliminated based on optimality conditions by examining the matrix of objective coefficients.

After applying the size reduction procedure, out of 1200 Sel-Col instances (to be described in the next section), 1131 became size-wise suitable for DA.
Table \ref{tab:selcol_reduction_stats} gives the percentage of instances that fit into DA before and after the size reduction operation, as well as the percentage reduction in the number of variables, separately for PG and ER instances.
We see that our heuristic procedure increases the number of instances eligible for DA by more than 60\% and 53\% for the two instance groups, by reducing their sizes around 85\% on the average.


\input{tables/SelCol_reduction_stats}

\subsubsection{Instances}

A Sel-Col instance is comprised of a graph and a partition of its vertex set, i.e., a set of clusters. 
Our instances are based on 400 distinct graph instances, with 200 being perfect graphs, and 200 \cite{erdos1959random} type graphs. 
The graphs contain 50 to 500 vertices and have average edge densities ranging from 0.1 to 0.7. 
The full set of Sel-Col instances are obtained by combining each graph with three different sets of clusters, where the number of vertices in each cluster vary between 2--5, 4--7, and 6--9. 
So, from each of the 400 graph instances, three different Sel-Col instances are generated, which makes a total of 1200 instances, where the ones based on perfect graphs are those used in \citep{cseker2020exact}.
Brief descriptive information about our test bed, which we obtained after applying the two-phase heuristic procedure explained above, is provided in Table \ref{tab:selcol_inst_info}.
The two rows of the table are reserved for the two groups of instances, PG and ER, and the columns respectively provide the number of instances, the range of the number of vertices, the edge densities of graphs, the range of the number of clusters, the range of the number of available colors in the IP and QUBO models (``\# avail colors"), and the range of the number of variables in these instances.

\input{tables/SelCol_instance_info}

\subsubsection{Experimental results}

We now present the results of our computational study that we conducted on our Sel-Col test bed.
As before, we employed DA, GUROBI, CPLEX and SCIP in our experiments, and used 30- and 60-second time limits.
We additionally tested GUROBI and CPLEX under a 10-minute limit to observe how much more they could achieve if longer solution times were allowed.
When solving the IP models with GUROBI and CPLEX, we have incorporated a set of symmetry breaking constraints, which are used in \citep{cseker2019decomposition,cseker2020exact}, in order to reduce the number of different representations of equivalent solutions in the feasible solution space. 
In addition, we fed the solution found through our two-phase heuristic as a warm-start to the IP solvers. 
In the figures that follow, we present the results of our experiments under 30-second time limit, but we note that the general picture does not change much in the 60-second experiments.
The points on the $x$-axes of figures represent individual instances sorted in ascending order with respect to the number of variables they contain.

\begin{figure}[ht]
    \centering
	\begin{subfigure}{0.49\textwidth}
		\centering
	    \includegraphics[scale=0.62]{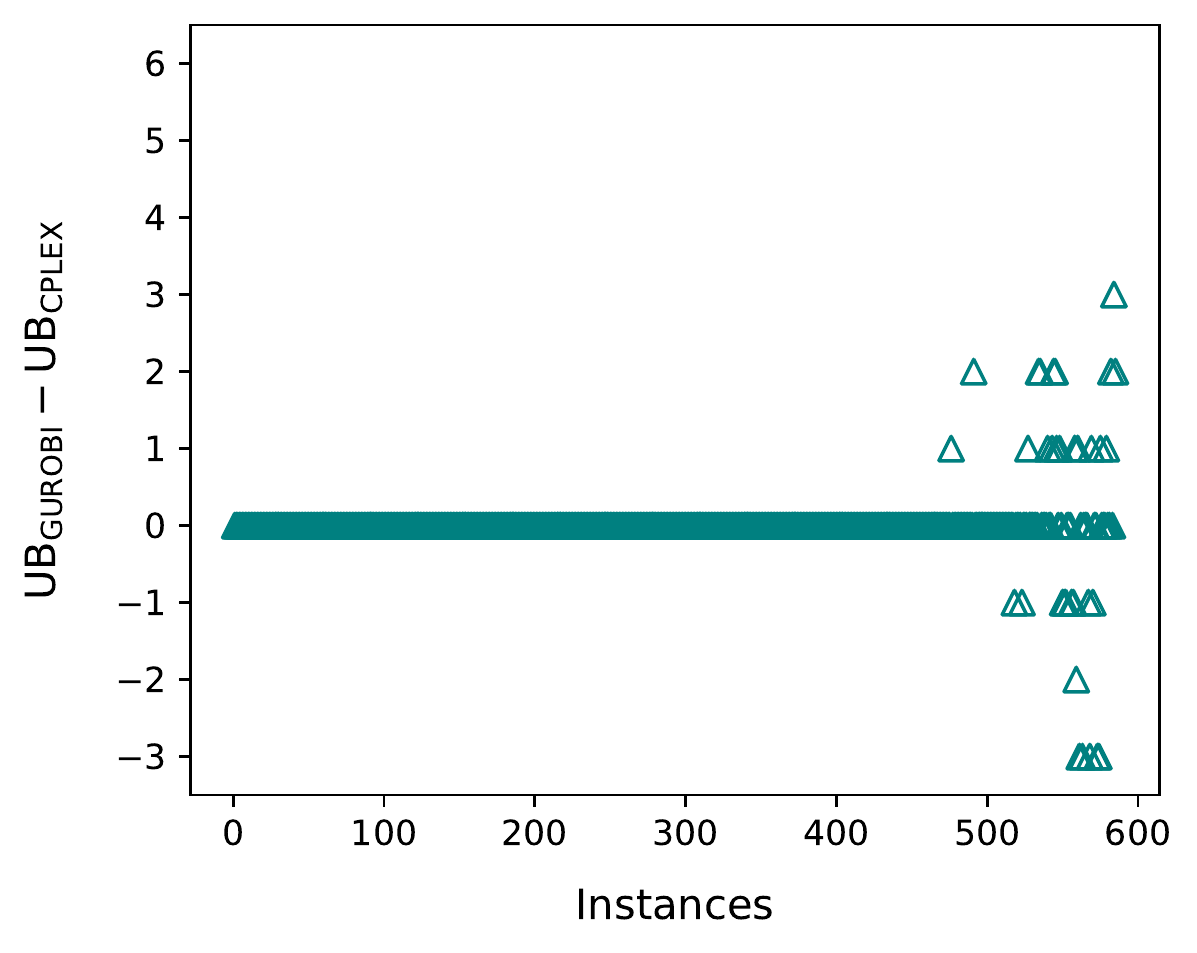}
		\caption{PG}
		\label{subfig:}
	\end{subfigure}
	\begin{subfigure}{0.49\textwidth}
		\centering
        \includegraphics[scale=0.62]{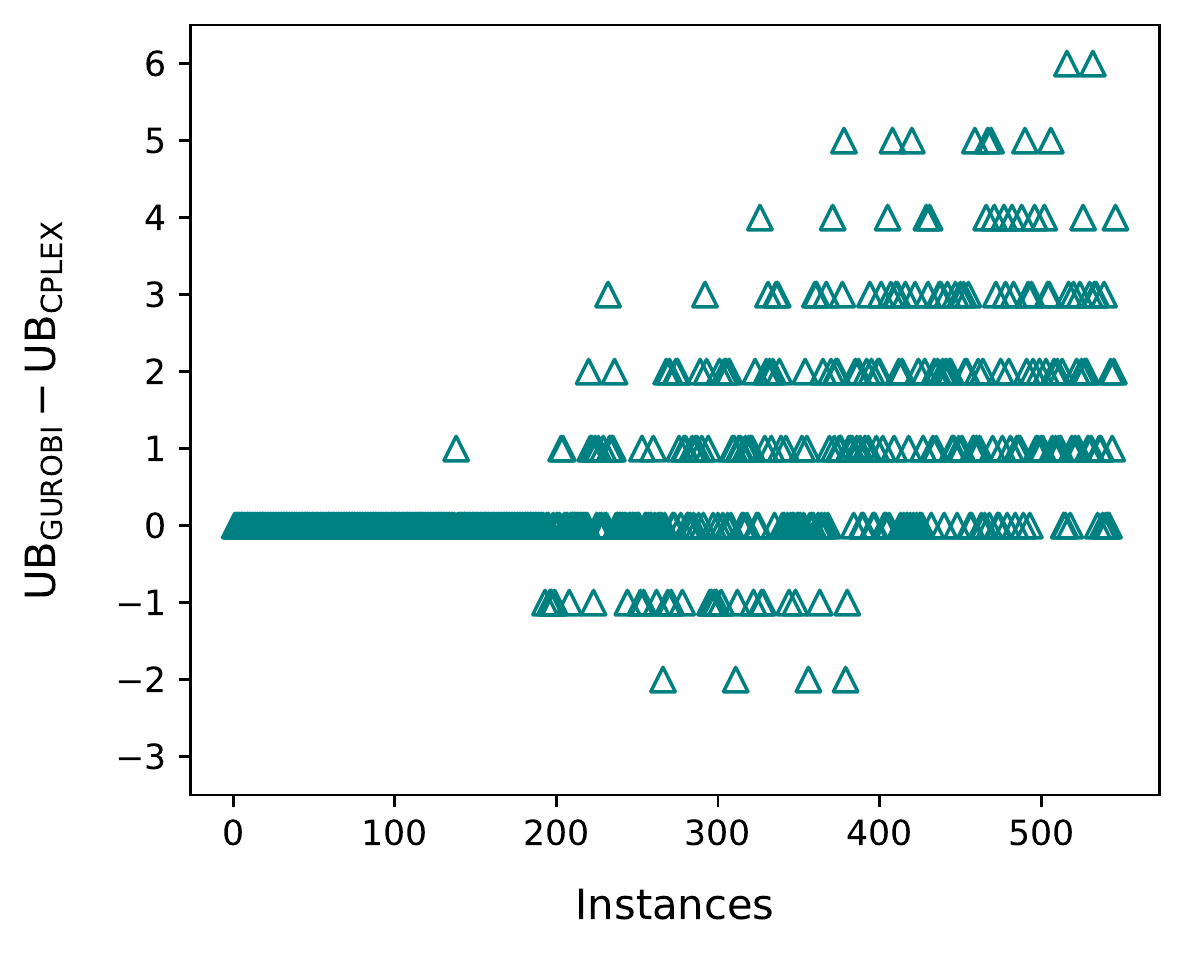}
        \caption{ER}
        \label{subfig:}
	\end{subfigure}
	\caption{Performance comparison of CPLEX and GUROBI for Sel-Col under 30-second time limit in terms of objective values (UBs).}
	\label{fig:selcol_gurobi_vs_cplex_UBdiff}
\end{figure}

{\partitlestyle Comparison of exact solvers.} 
We start with comparing the performances of GUROBI and CPLEX on the basis of the upper bounds they provide, i.e., the objective values of the best feasible solutions found.
We note that this measure is interesting for this problem class since the optimal values (i.e., the minimum number of colors needed) are usually small, and it thus provides a better idea on the quality of the obtained solutions than the optimality gaps.
Figure \ref{fig:selcol_gurobi_vs_cplex_UBdiff} presents the difference between the upper bounds from the two solvers.
The number of instances being high, we plot the differences rather than the individual upper bounds in order to make the exposition of the comparison easily interpretable.
Since we subtract the upper bound values from CPLEX from those of GUROBI, a data point above the zero level means that GUROBI is outperformed by CPLEX for the associated instance. 
Conversely, a data point lying below the zero level indicates that GUROBI outperforms for that instance. 

For PG instances, GUROBI and CPLEX yield the same objective values in about 94\% of the instances (zero level), and there is no obvious outperformance of either in the rest.
For ER instances, on the other hand, CPLEX shows a better performance than GUROBI by yielding smaller and same objective values in nearly 37\% and 58\% of the instances, respectively.

\begin{figure}[ht]
    \centering
	\begin{subfigure}{0.49\textwidth}
		\centering
	    \includegraphics[scale=0.62]{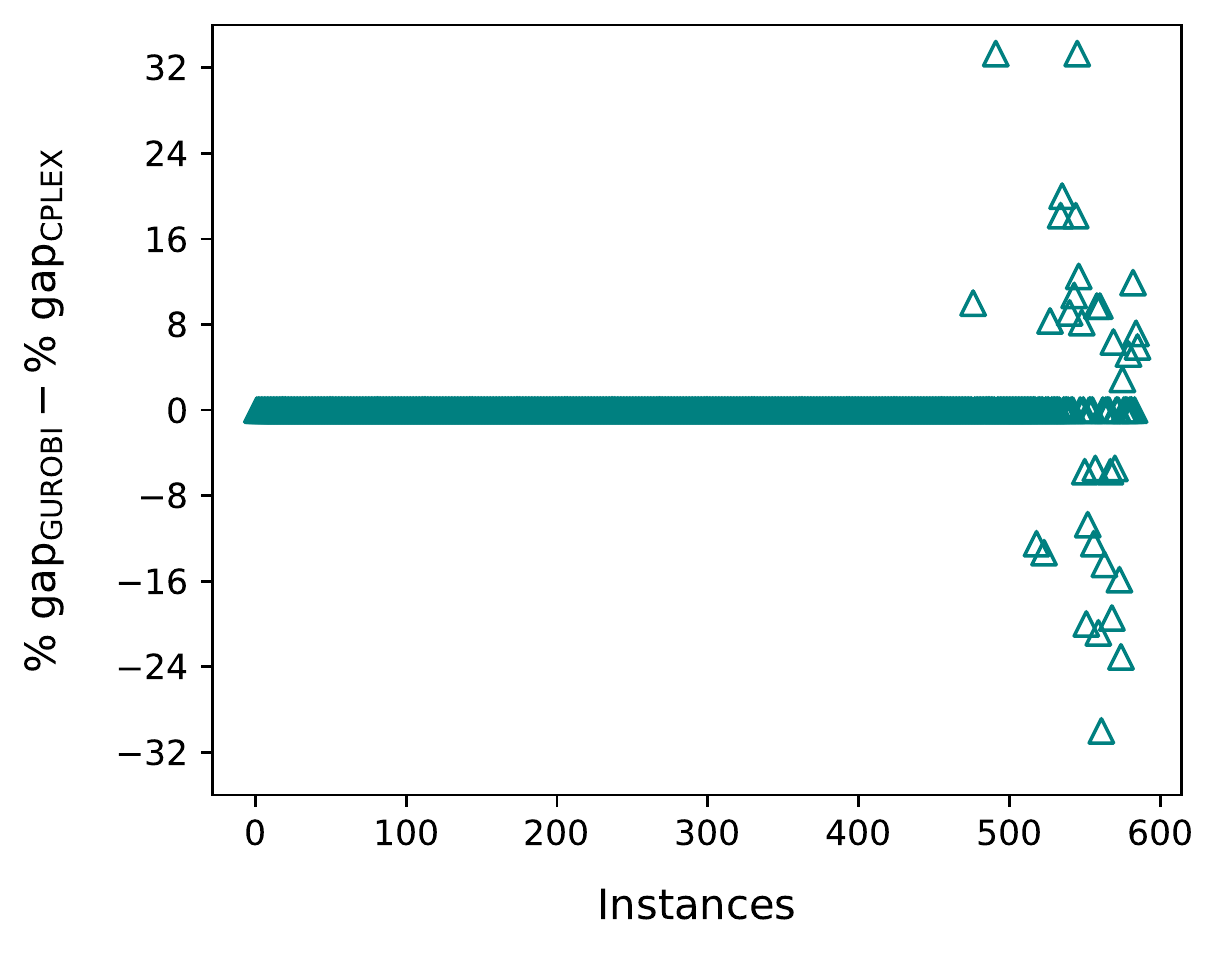}
		\caption{PG}
		\label{subfig:}
	\end{subfigure}
	\begin{subfigure}{0.49\textwidth}
		\centering
        \includegraphics[scale=0.62]{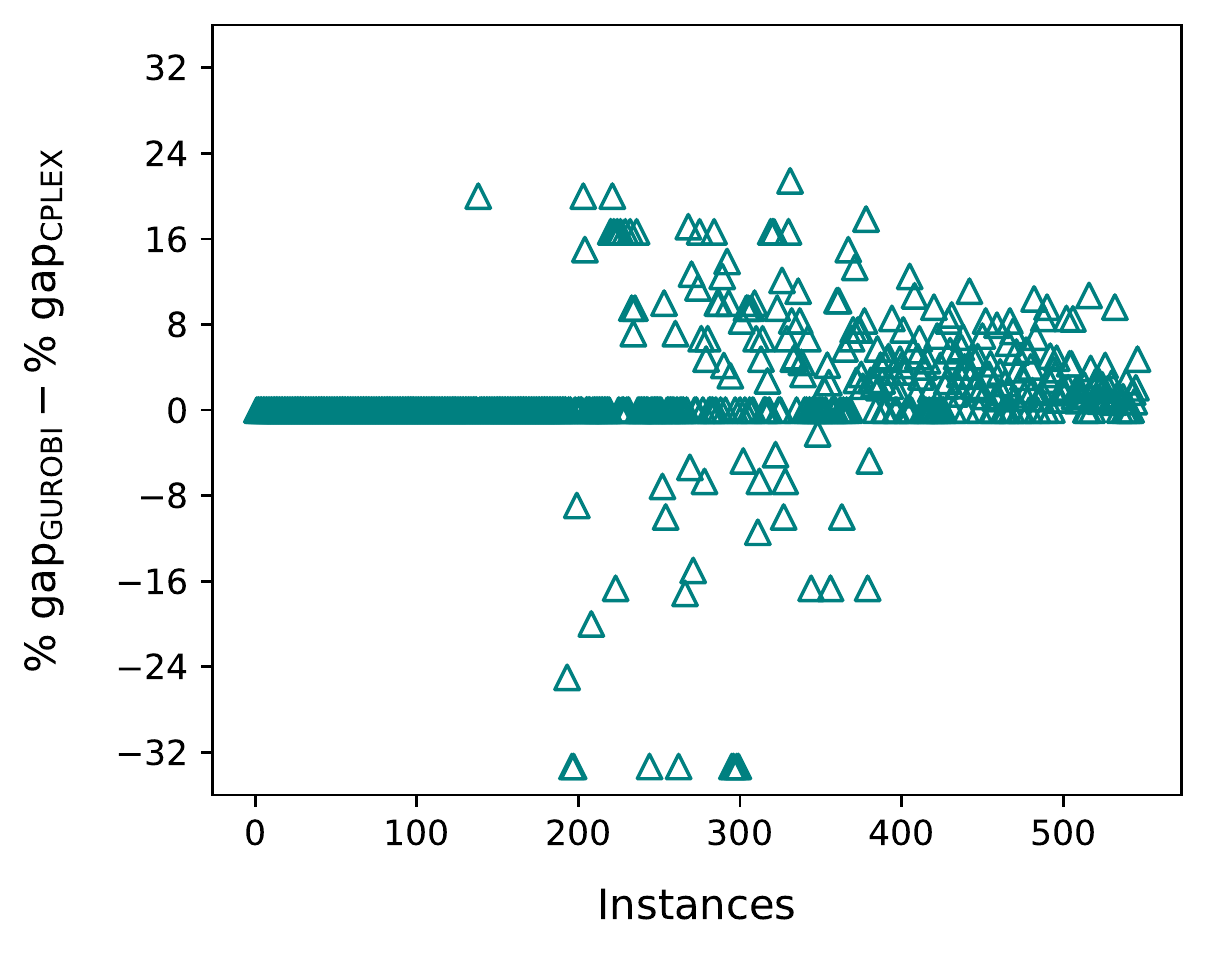}
        \caption{ER}
        \label{subfig:}
	\end{subfigure}
	\caption{Performance comparison of CPLEX and GUROBI for Sel-Col under 30-second time limit in terms of percentage gaps with respect to the best lower bound values.}
	\label{fig:selcol_gurobi_vs_cplex_ownGap}
\end{figure}

The upper bound differences provide a comparison in terms of the quality of the feasible solutions provided, but do not say much on how close the solutions are to the optimals, because the lower bounds are absent in the analysis. 
Therefore, we next compare the performances using the difference in percentage optimality gaps with respect to the best lower bound values available (``\% gap") in Figure \ref{fig:selcol_gurobi_vs_cplex_ownGap}. 
As before, we plot the differences, rather than the individual percentage gap values separately for each solver, in order to avoid a convoluted picture arising from overlapping data points. 
Our observation is that the percentage gap values of GUROBI and CPLEX are the same in the majority of instances (zero level), and in most of the remaining, CPLEX yields better results (above the zero level). 
This indicates that CPLEX usually yields better solutions than GUROBI, and the difference is more evident in ER instances. 
Hence, we continue our evaluation using the results of CPLEX. 

\begin{figure}[h]
    \centering
	\begin{subfigure}{0.49\textwidth}
		\centering
	    \includegraphics[scale=0.62]{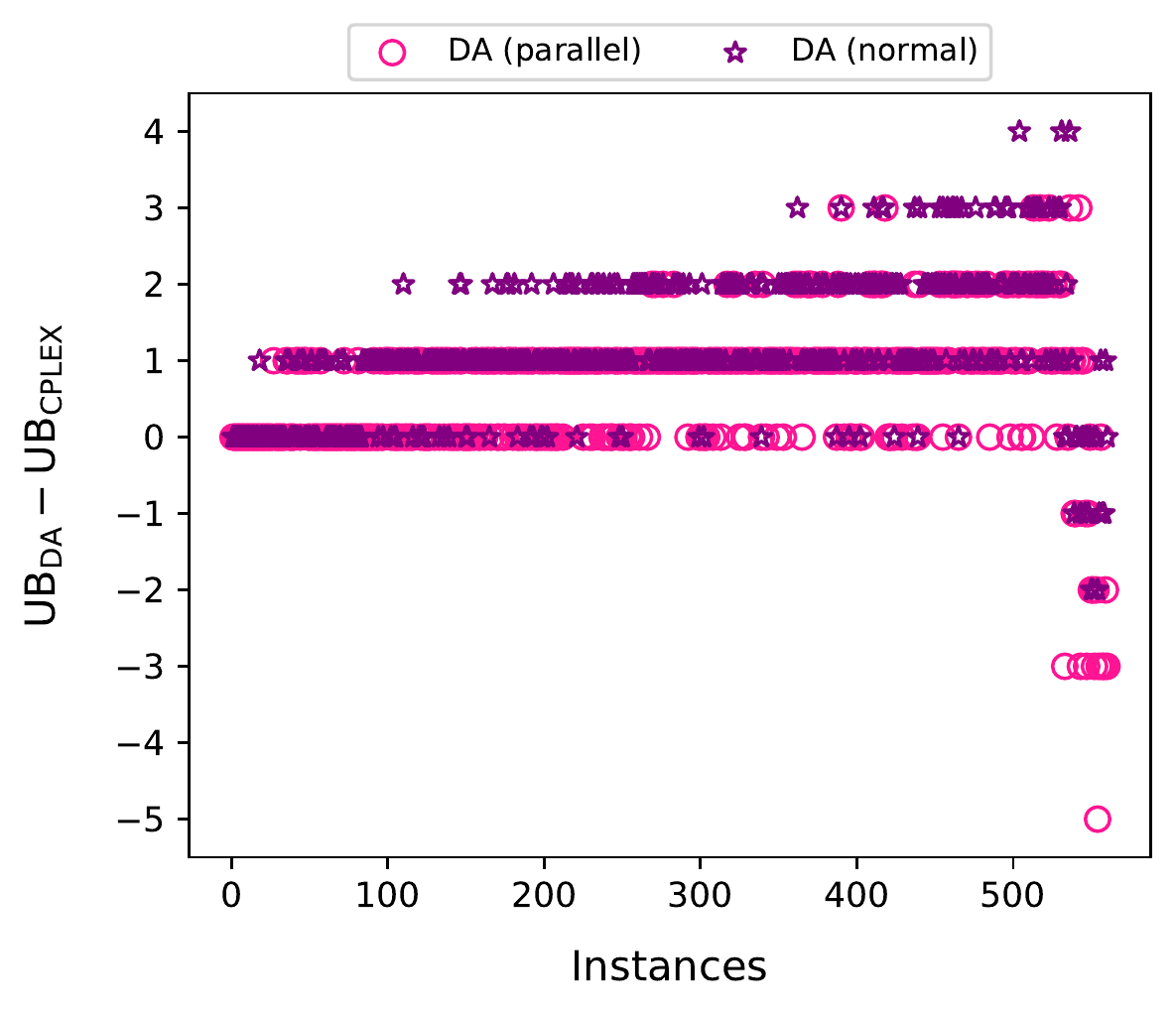}
		\caption{PG}
		\label{subfig:}
	\end{subfigure}
	\begin{subfigure}{0.49\textwidth}
		\centering
        \includegraphics[scale=0.62]{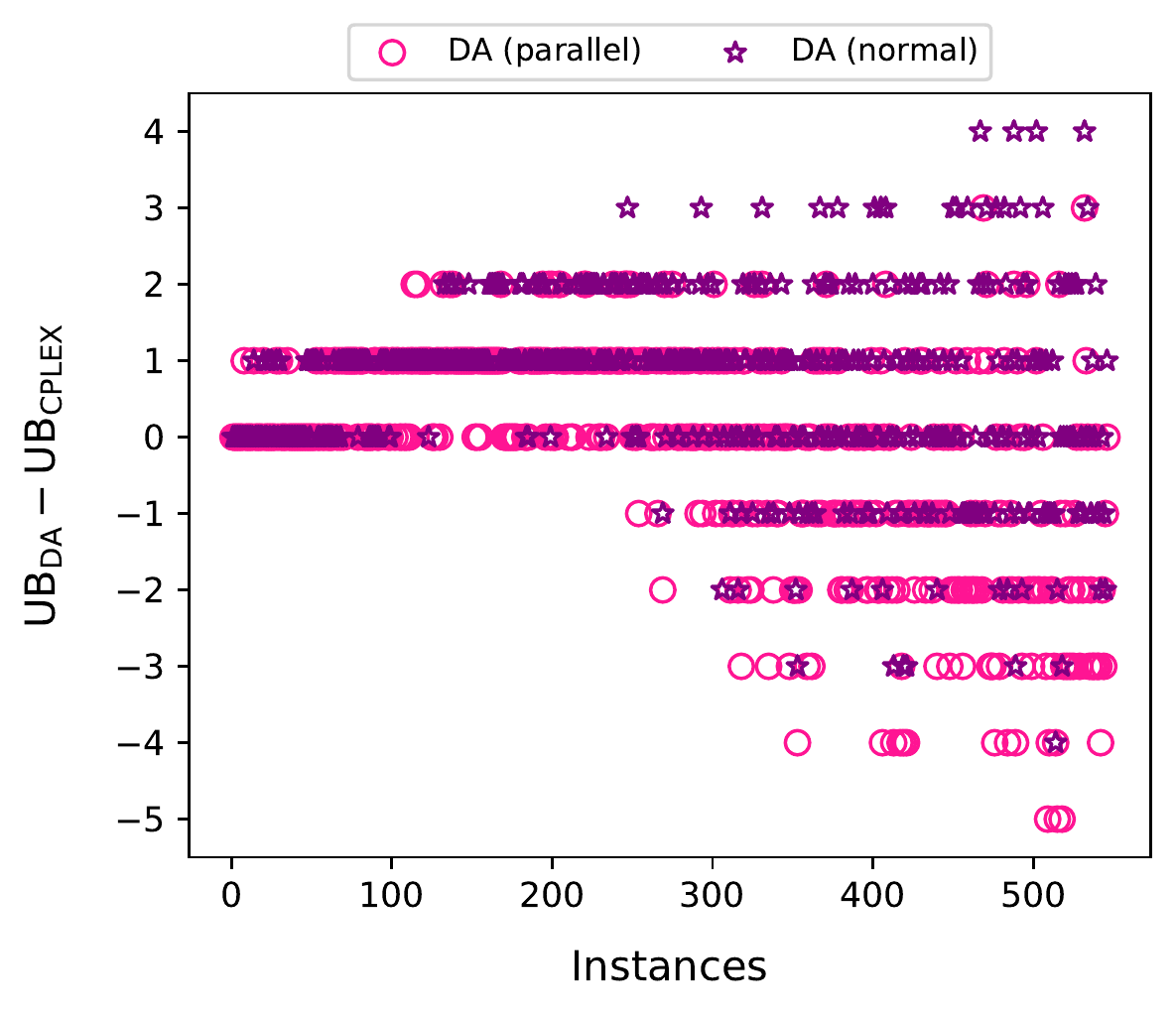}
        \caption{ER}
        \label{subfig:}
	\end{subfigure}
	\caption{Performance comparison of DA and CPLEX for Sel-Col under 30-second time limit.}
	\label{fig:selcol_da2_vs_cplex_UBdiff}
\end{figure}

{\partitlestyle Comparison of DA and CPLEX.} 
We now compare DA's performance to CPLEX. 
Figure \ref{fig:selcol_da2_vs_cplex_UBdiff} presents the difference of the DA's objective values in normal and parallel modes from those of the CPLEX.

For the majority of PG instances, the quality of the objective values found by DA are at most two units worse than those of CPLEX, but the overall performance of CPLEX is better than both modes of DA.
In particular, normal and parallel modes of DA yield at least as good results in nearly 21\% and 39\% of the PG instances.
For ER instances, however, normal mode of DA performs at least as well as CPLEX in 39\% of the instances, while it becomes more than 61\% for the parallel mode.
So, parallel mode of DA is successful in yielding good-quality solutions especially for ER instances.

\begin{figure}[h]
    \centering
	\begin{subfigure}{0.49\textwidth}
		\centering
	    \includegraphics[scale=0.62]{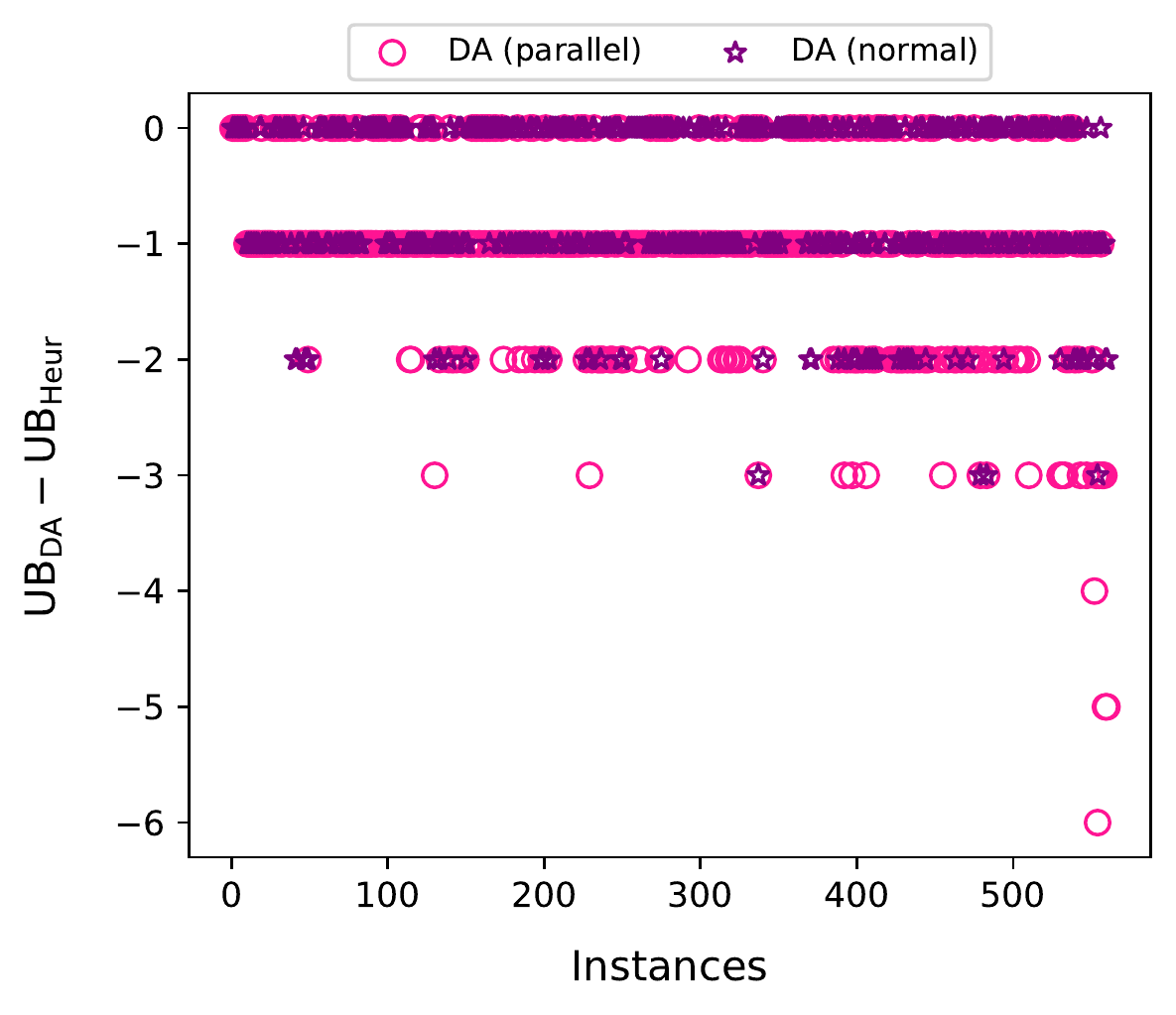}
		\caption{PG}
		\label{subfig:}
	\end{subfigure}
	\begin{subfigure}{0.49\textwidth}
		\centering
        \includegraphics[scale=0.62]{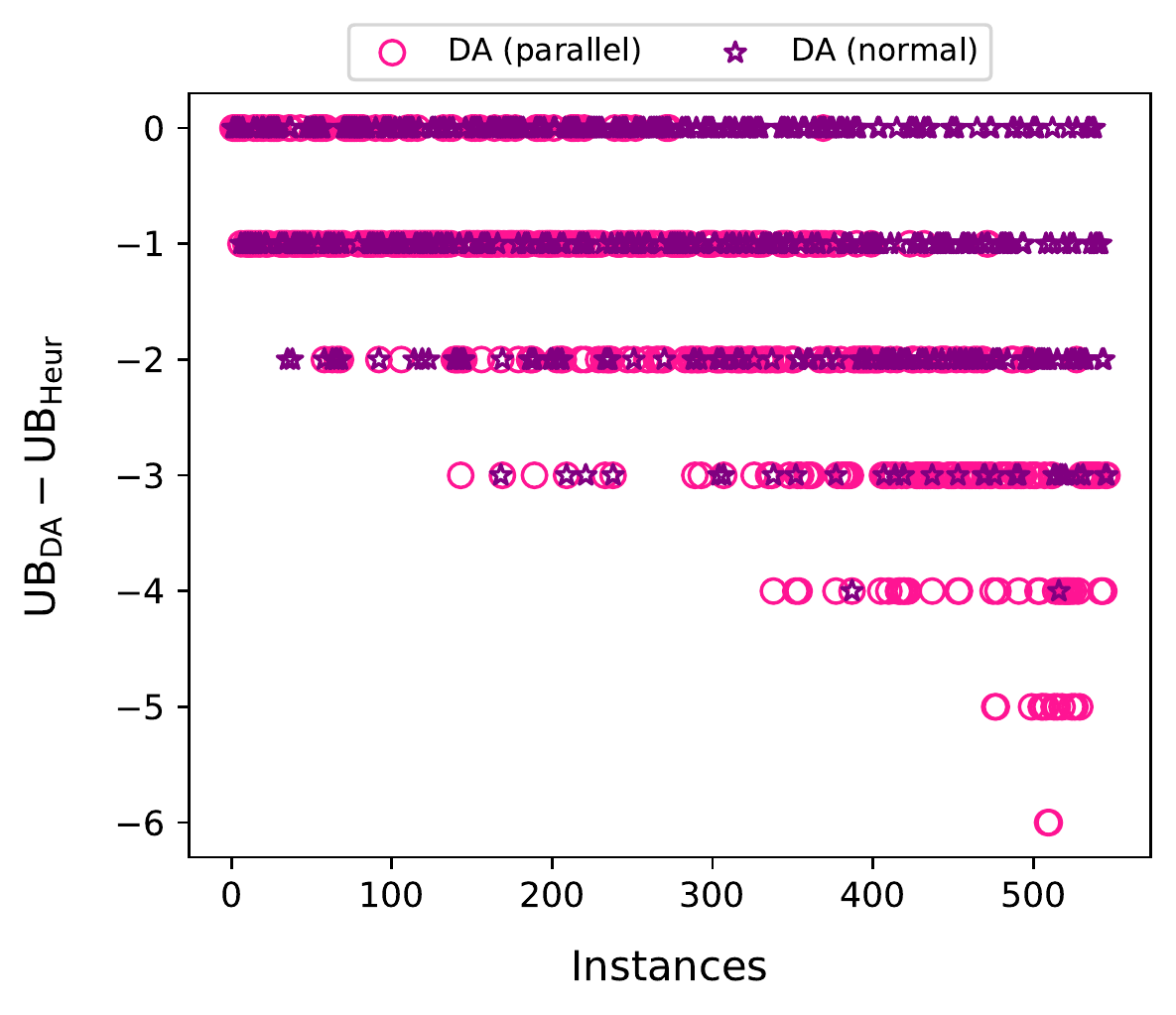}
        \caption{ER}
        \label{subfig:}
	\end{subfigure}
	\caption{Marginal contribution of DA (under 30-second time limit) upon the values two-phase heuristic achieves for Sel-Col.}
	\label{fig:selcol_da2_vs_heur_UBdiff}
\end{figure}

Having evaluated the performance of DA in comparison to CPLEX, we now investigate the marginal contribution of DA on top of the objective values achieved by our two-phase heuristic.
Figure \ref{fig:selcol_da2_vs_heur_UBdiff} shows the difference of DA's objective values from those of the two-phase heuristic.

We observe that the amount of improvement, if any, ranges from one to six units, and the parallel mode of DA performs somewhat better than the normal mode.
Specifically, DA improves upon the two-phase heuristic in nearly 50\% and 78\% of PG instances with normal and parallel modes, respectively, while these percentages increase to 62\% and 88\% for ER instances.
So, especially the parallel mode of DA improves the solution quality in the majority of instances, if there is any room for it, and the improvement manifests more evidently in ER instances.

\begin{figure}[ht]
    \centering
	\begin{subfigure}{0.49\textwidth}
		\centering
	    \includegraphics[scale=0.62]{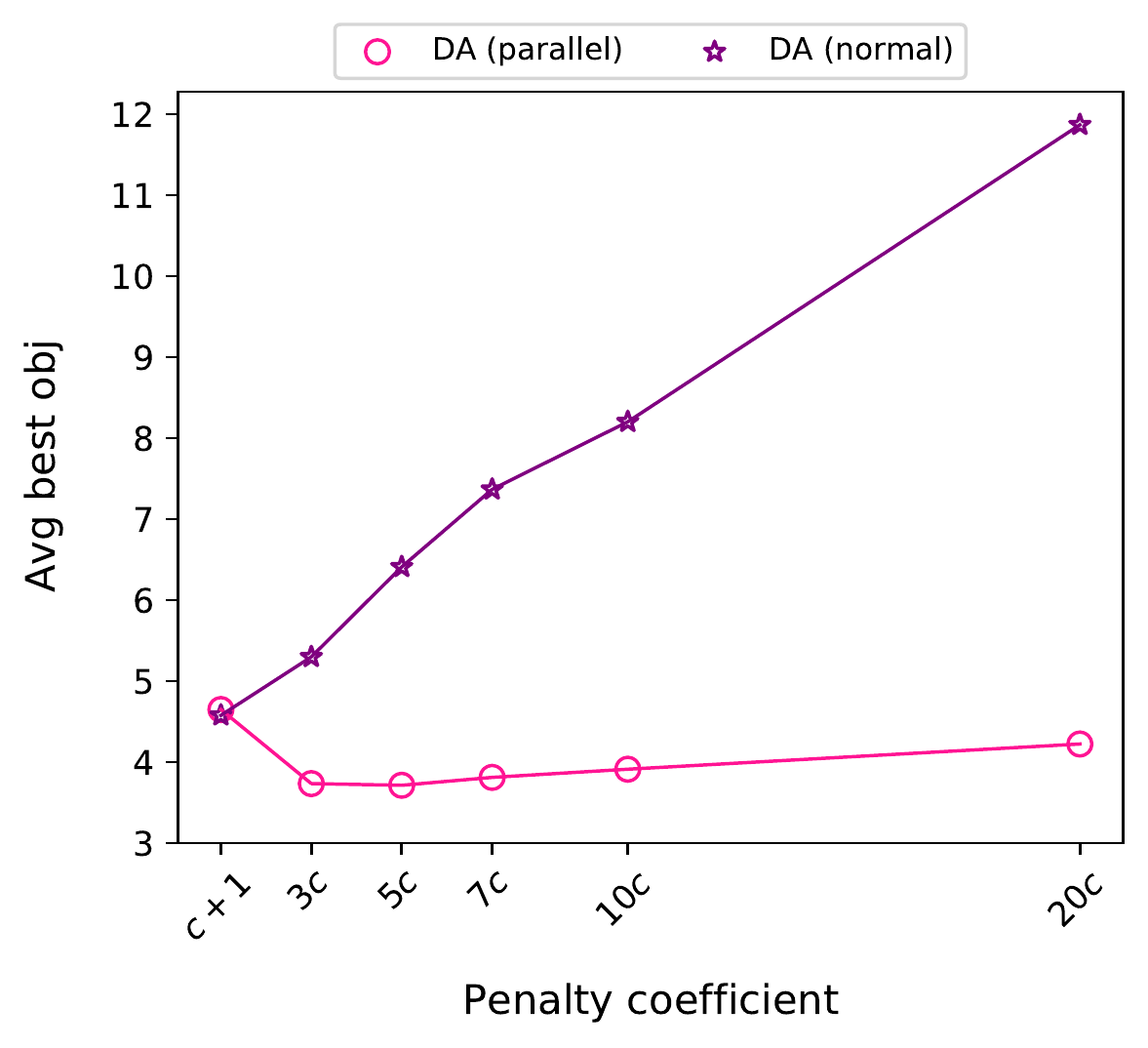}
		\caption{PG}
		\label{subfig:selcol_pen_trials_pg}
	\end{subfigure}
	\begin{subfigure}{0.49\textwidth}
		\centering
        \includegraphics[scale=0.62]{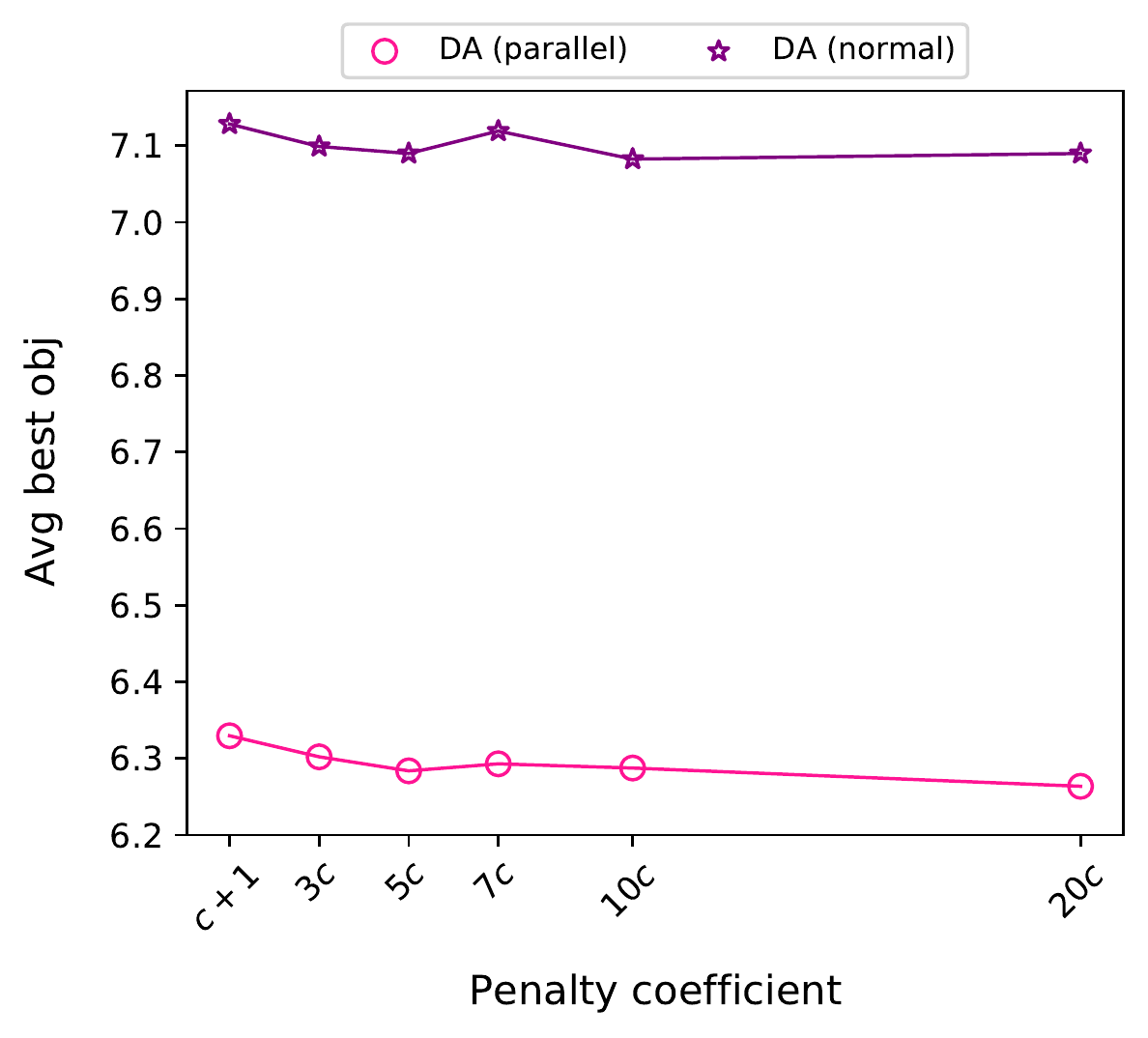}
        \caption{ER}
        \label{subfig:selcol_pen_trials_er}
	\end{subfigure}
	\caption{Penalty coefficient analysis for Sel-Col under 30-second time limit.}
	\label{fig:selcol_pen_trials}
\end{figure}

{\partitlestyle Penalty coefficient analysis.}
Next, we investigate the effect of the penalty coefficient values on DA's performance. 
For Sel-Col, the maximum value that the objective function can take is the number of available colors. 
Moreover, if a solution to the QUBO model in \eqref{m:selsol_qubo} violates some of the original constraints in \eqref{eq:selcol_IP_c1}--\eqref{eq:selcol_IP_c3}, the magnitude of violation is at least one. 
Then, a penalty coefficient ($\penaltyCoef$) value strictly larger than the number of available colors leads to a penalty greater than the largest possible objective value any feasible QUBO solution can yield, and hence the value of the QUBO objective for an optimal solution always lies below those that are infeasible for the IP model, rendering the QUBO formulation in \eqref{m:selsol_qubo} exact.
Let $\numColors$ denote the number of available colors (which can at most be equal to the number of clusters $P$).
Though $\penaltyCoef = \numColors +1$ is sufficient to make the QUBO model exact, we experimented with various penalty coefficient values as a function of the number of available colors $\numColors$ to test the effect on the performance of DA.
Figure \ref{fig:selcol_pen_trials} shows the average objective values from both modes of DA for different penalty coefficient settings.

We see from Figures \ref{subfig:selcol_pen_trials_pg} and \ref{subfig:selcol_pen_trials_er} that the effect of the penalty coefficient value manifests differently in normal and parallel modes, and in PG and ER instances.
Even though smaller penalty coefficient values typically lead to improved performance in DA, as analyzed and reported in \citep{cseker2020routing} and mentioned in \citep{cohen2020constrainedclust,cohen2020ising}, it appears from our results in Sel-Col context that higher penalty coefficient values may ameliorate or deteriorate performance, and no strong relationship can be inferred.
In our experiments, we used $\penaltyCoef = 5\numColors$, which we observed to be the best overall option to get the best out of DA.
We note that even if we did not make this sensitivity analysis for DA and simply presented the results with $\penaltyCoef = \numColors + 1$, our overall conclusions would not change much, because the marginal improvement it creates for DA does not dominate the performance gap between DA and other solvers.

{\partitlestyle Overall performances.} 
We conclude this section with tabular summaries of the overall results.
Table \ref{tab:selcol_summary} reports the averages of different performance measures over the PG and ER instances under different time limits.
There are two sets of rows, reserved for the results associated with PG and ER instances.
For each group of instances, 
the average of upper bounds (``UB"), 
percentage of solutions that are feasible (``\% feas"),
average solution time (``Time (sec)"), 
average percentage gap with respect to the best lower bound available (``\% gap"), and
average percentage gap of the solver (``\% solver gap"), i.e., gap values using the lower bound values from that solver (rather than the best values) are given.
There are four groups of columns, one dedicated for the normal mode of DA, and the other three reserved for the results obtained under the 30-second, 60-second and 10-minute time limits.
These column groups contain the results from solving the QUBO formulation with DA and/or solving the IP formulations with GUROBI, CPLEX and SCIP.  
The reason for the normal mode of DA to comprise a separate group is that no matter how high the number iterations are set, it takes less than a second to solve an instance. 

\input{tables/SelCol_summary_IP}

Our observation from Table \ref{tab:selcol_summary} is that under both 30-second and 60-second time limits, the parallel mode of DA yields the best upper bounds for ER instances, and CPLEX is only slightly worse in that regard. 
Moreover, the upper bounds from 10-minute experiments with CPLEX and GUROBI are not much smaller than those of DA, especially for ER instances.  
In terms of solution times, GUROBI and CPLEX yield better results than the parallel mode of DA {\cfujitsu run with default hyperparameters} (``DA (p)").
However, the normal mode of DA (``DA (n)") is extremely fast, and is able to yield good-quality results in less than half a second.
For ER instances in particular, we see that the upper bounds from the normal mode of DA are only about one unit away from the 10-minute values of GUROBI and CPLEX.
Noting that DA always yields feasible solutions except for a few PG instances, which we see from the ``\% feas" rows, we can say that the normal mode of DA is by far the best option to obtain good-quality solutions in a very short amount of time, such that GUROBI and CPLEX typically fail to complete even presolve phase of the solution process within that time frame.
As for the parallel mode of DA, we see that it is competitive with GUROBI or CPLEX especially in ER instances, which is where they struggle more.

Lastly, we present the results of our experiments where we utilized the three exact solvers to solve the QUBO models for Sel-Col.
Table \ref{tab:selcol_summary_qubo} contains the average upper bound, percentage feasibility of solutions, solution time, and percentage gap values from solvers and with respect to the best available lower bound values for GUROBI, CPLEX, and SCIP under 30-second and 60-second time limits.
The general structure of this table is the same as that of Table \ref{tab:selcol_summary}.
The \largeVal~signs stand for values greater than $10^5$.

\input{tables/SelCol_summary_QUBO}

We see from the ``\% feas" rows of the table that the feasibility percentages vary between 11\% and 85\%.
Even though the feasibility percentages are not satisfactory in general, the ones from SCIP are far from being at an acceptable level.
Since infeasibilities are penalized in the objective function and the penalty coefficient values we use for Sel-Col ensure that the objective values of infeasible solutions are always larger than those of feasible solutions, the average upper bound values are much larger than what they are expected to be even in the worst case. 
Indeed, the poor performance GUROBI and CPLEX as QUBO solvers is expected, because these solvers are specialized in solving linear and integer linear programs particularly, which is confirmed by their performance in solving the IP formulations, presented above.
Overall, we can say that as QUBO solvers, GUROBI and CPLEX are the best among the three alternatives, and the performance SCIP of lies far below the others, which is particularly evident from the percentage of feasible solutions.
We note here that even though some of these solvers occasionally yielded invalid lower bounds for quadratic models, which we have observed for some problem classes as mentioned before, we did not obtain any inconsistent results in our Sel-Col experiments.
Though this does not prove the validity of the percentage gap values, we nevertheless report them to make our summary table complete and give a general idea on the performances with respect to those measures.

%% file: figures/SelCol/illustrativeExample.tex
\begin{figure}[ht]
	\centering
	\scalebox{0.6}[0.6]{ 
		\begin{tikzpicture}[main_node/.style={circle,fill=white!80,draw,inner sep=0pt, minimum size=18pt}, line width=0.75pt]	
		\node[main_node] (v1) at (-4,2) {1};
		\node[main_node] (v2) at (-1,2) {2};
		\node[main_node] (v3) at (-1,-1) {3};
		\node[main_node] (v4) at (-4,-1) {4};
		\node[main_node] (v5) at (-3.25,1.25) {5};
		\node[main_node] (v6) at (-1.75,1.25) {6};
		\node[main_node] (v7) at (-1.75,-0.25) {7};
		\node[main_node] (v8) at (-3.25,-0.25) {8};
		\draw (v1) -- (v2) -- (v3) -- (v4) -- (v1) -- (v5) -- (v6) -- (v2);
		\draw (v6) -- (v7) -- (v3);
		\draw (v7) -- (v8) -- (v4);
		\draw (v5) -- (v8);
		\draw[dashed, draw=red, rotate around={38:(-3.75,1.75)}] (-3.75,1.75) ellipse (0.55 and 1.25);
		\draw[dashed, draw=red, rotate around={-38:(-1.25,1.75)}]  (-1.25,1.75) ellipse (0.55 and 1.25);
		\draw[dashed, draw=red, rotate around={38:(-1.25,-0.75)}]  (-1.25,-0.75) ellipse (0.55 and 1.25);
		\draw[dashed, draw=red, rotate around={-38:(-3.75,-0.75)}]  (-3.75,-0.75) ellipse (0.55 and 1.25);
		\node at (-5,2.5) {\textcolor{red}{$V_1$}};
		\node at (0,2.5) {\textcolor{red}{$V_2$}};
		\node at (-5,-1.5) {\textcolor{red}{$V_3$}};
		\node at (0,-1.5) {\textcolor{red}{$V_4$}};
		\node[main_node,fill=gray!40] (v9) at (4,2) {1};
		\node[main_node,fill=gray!100] (v10) at (7,2) {2};
		\node[main_node,fill=gray!40] (v11) at (7,-1) {3};
		\node[main_node,fill=gray!100] (v12) at (4,-1) {4};
		\node[main_node] (v13) at (4.75,1.25) {5};
		\node[main_node] (v14) at (6.25,1.25) {6};
		\node[main_node] (v15) at (6.25,-0.25) {7};
		\node[main_node] (v16) at (4.75,-0.25) {8};
		\draw (v9) -- (v10) -- (v11) -- (v12) -- (v9) -- (v13) -- (v14) -- (v10);
		\draw (v14) -- (v15) -- (v11);
		\draw (v15) -- (v16) -- (v12);
		\draw (v13) -- (v16);
		\draw[dashed, draw=red, rotate around={38:(4.25,1.75)}] (4.25,1.75) ellipse (0.55 and 1.25);
		\draw[dashed, draw=red, rotate around={-38:(6.75,1.75)}]  (6.75,1.75) ellipse (0.55 and 1.25);
		\draw[dashed, draw=red, rotate around={38:(6.75,-0.75)}]  (6.75,-0.75) ellipse (0.55 and 1.25);
		\draw[dashed, draw=red, rotate around={-38:(4.25,-0.75)}]  (4.25,-0.75) ellipse (0.55 and 1.25);
		\node at (3,2.5) {\textcolor{red}{$V_1$}};
		\node at (8,2.5) {\textcolor{red}{$V_2$}};
		\node at (3,-1.5) {\textcolor{red}{$V_3$}};
		\node at (8,-1.5) {\textcolor{red}{$V_4$}};
		\node[main_node,fill=gray!100] (v17) at (12,2) {1};
		\node[main_node] (v18) at (15,2) {2};
		\node[main_node,fill=gray!100] (v19) at (15,-1) {3};
		\node[main_node] (v20) at (12,-1) {4};
		\node[main_node] (v21) at (12.75,1.25) {5};
		\node[main_node,fill=gray!100] (v22) at (14.25,1.25) {6};
		\node[main_node] (v23) at (14.25,-0.25) {7};
		\node[main_node,fill=gray!100] (v24) at (12.75,-0.25) {8};
		\draw (v17) -- (v18) -- (v19) -- (v20) -- (v17) -- (v21) -- (v22) -- (v18);
		\draw (v22) -- (v23) -- (v19);
		\draw (v23) -- (v24) -- (v20);
		\draw (v21) -- (v24);
		\draw[dashed, draw=red, rotate around={38:(12.25,1.75)}] (12.25,1.75) ellipse (0.55 and 1.25);
		\draw[dashed, draw=red, rotate around={-38:(14.75,1.75)}]  (14.75,1.75) ellipse (0.55 and 1.25);
		\draw[dashed, draw=red, rotate around={38:(14.75,-0.75)}]  (14.75,-0.75) ellipse (0.55 and 1.25);
		\draw[dashed, draw=red, rotate around={-38:(12.25,-0.75)}]  (12.25,-0.75) ellipse (0.55 and 1.25);
		\node at (11,2.5) {\textcolor{red}{$V_1$}};
		\node at (16,2.5) {\textcolor{red}{$V_2$}};
		\node at (11,-1.5) {\textcolor{red}{$V_3$}};
		\node at (16,-1.5) {\textcolor{red}{$V_4$}};
		\node at (-2.5,-2.75) {\large{(a)}};
		\node at (5.5,-2.75) {\large{(b)}};	
		\node at (13.5,-2.75) {\large{(c)}};
		%
		\node at (-2.5,3) {};
		\end{tikzpicture}
	}
	\caption{(a) A graph $G$ with a partition of its vertex set into four clusters shown in dashed ellipses, (b) an optimally colored selection $\{1,2,3,4\}$ in $G$, (c) an optimal selection $\{1,6,3,8\}$ in $G$ with an optimal coloring of it (taken from \citep{cseker2020exact} with permission).	\label{fig:selcol_example}}
	{}
\end{figure}

%% file: tables/SelCol_reduction_stats.tex
\begin{table}[htbp]
  \centering
  \caption{Size reduction statistics for Sel-Col instances.}
  \def\arraystretch{1}
  \scalebox{0.86}{
    \begin{tabular}{cccc}
    \toprule
      & \multicolumn{2}{c}{\% eligible for DA} & \parbox[t]{2cm}{\centering \%} \\
      \cmidrule(lr){2-3}
      \parbox[t]{1.2cm}{\centering \phantom{o}} & \parbox[t]{1.5cm}{\centering Before} & \parbox[t]{1.5cm}{\centering After} & \parbox[t]{2cm}{\centering reduction} \\
    \midrule
    PG    & 37.0 & 97.5 & 88.8 \\[0.05cm]
    ER    & 37.2 & 91.0 & 80.6 \\[0.05cm]
    \bottomrule
    \end{tabular}%
    }
  \label{tab:selcol_reduction_stats}%
\end{table}%

%% file: tables/SelCol_instance_info.tex
\begin{table}[htbp]
  \centering
  \caption{Sel-Col instance information.}
  \def\arraystretch{1}
  \scalebox{0.85}{
    \begin{tabular}{ccccccc}
    \toprule
      \parbox[t]{1cm}{\centering \phantom{oo}} & \parbox[t]{2.1cm}{\centering \# instances} & \parbox[t]{2.1cm}{\centering \# vertices} & \parbox[t]{2.4cm}{\centering Density} & \parbox[t]{2.1cm}{\centering \# clusters} & \parbox[t]{2.5cm}{\centering \# avail colors} & \parbox[t]{2.1cm}{\centering \# variables} \\
    \midrule
    PG    & 586   & 50--500 & 0.1, 0.3, 0.5, 0.7 & 6--150 & 1--22 & 51--7722 \\[0.05cm]
    ER    & 547   & 50--500 & 0.1, 0.3, 0.5, 0.7 & 6--152 & 1--25 & 51--8118 \\
    \bottomrule
    \end{tabular}%
    }
  \label{tab:selcol_inst_info}%
\end{table}%

%% file: tables/SelCol_summary_IP.tex
\begin{table}[htbp]
  \centering
    \caption{Average performance of DA, GUROBI, CPLEX and SCIP as IP solvers for Sel-Col.}
  \def\arraystretch{1}
  \scalebox{0.79}{
    \begin{tabular}{cc rrrrrrrrr}
    \toprule
    & & & \multicolumn{3}{c}{30-sec limit} & \multicolumn{3}{c}{60-sec limit} & \multicolumn{2}{c}{10-min limit} \\
    \cmidrule(lr){4-6}\cmidrule(lr){7-9}\cmidrule(lr){10-11}
    \parbox[t]{0.75cm}{\centering \phantom{o} } 
    & \parbox[t]{1.5cm}{\centering \phantom{o} } 
    & \parbox[t]{1.3cm}{\centering DA (n)} 
    & \parbox[t]{1.3cm}{\centering DA (p)} 
    & \parbox[t]{1.4cm}{\centering GUROBI} 
    & \parbox[t]{1.4cm}{\centering CPLEX} 
    & \parbox[t]{1.3cm}{\centering DA (p)}  
    & \parbox[t]{1.4cm}{\centering GUROBI} 
    & \parbox[t]{1.4cm}{\centering CPLEX} 
    & \parbox[t]{1.4cm}{\centering GUROBI} 
    & \parbox[t]{1.4cm}{\centering CPLEX} \\
    \midrule
    \multirow{5}[0]{*}{PG} 
          & UB & 6.4  & 3.7  & 2.9  & 2.9  & 3.7  & 2.9  & 2.9  & 2.7  & 2.7 \\[0.05cm]
          & \% feas & 96.7  & 99.5  & 100.0 & 100.0 & 99.5  & 100.0 & 100.0 & 100.0 & 100.0 \\[0.05cm]
          & Time (sec) & 0.4  & 31.9 & 5.9  & 5.5  & 63.1 & 8.7  & 8.3  & 38.1 & 38.9 \\[0.05cm]
          & \% gap  & 25.6 & 25.7 & 1.9  & 1.9  & 25.6 & 1.3  & 1.4  & 0.4  & 0.4 \\[0.05cm]
          & \% solver gap  &  {--}  &  {--}  & 5.6  & 4.6  &  {--}  & 3.7  & 3.2  & 0.7  & 0.6 \\[0.05cm]
    \midrule
    \multirow{5}[0]{*}{ER} 
          & UB & 7.1  & 6.3  & 7.1  & 6.4  & 6.3  & 6.9  & 6.3  & 6.0  & 6.0 \\[0.05cm]
          & \% feas & 100.0  & 100.0  & 100.0 & 100.0 & 100.0  & 100.0 & 100.0 & 100.0 & 100.0 \\[0.05cm]
          & Time (sec) & 0.4  & 31.3 & 20.2 & 21.5 & 62.2 & 38.5 & 41.3 & 343.9 & 357.2 \\[0.05cm]
          & \% gap & 46.8 & 46.8 & 40.2 & 38.8 & 46.8 & 38.8 & 37.7 & 35.2 & 35.3 \\[0.05cm]
          & \% solver gap &  {--}  &  {--}  & 49.2 & 47.6 &  {--}  & 45.8 & 45.2 & 36.4 & 36.2 \\[0.05cm]
    \bottomrule
    \end{tabular}%
    }
  \label{tab:selcol_summary}%
\end{table}%

%% file: tables/SelCol_summary_QUBO.tex
\begin{table}[htbp]
  \centering
  \caption{Average performance of GUROBI, CPLEX and SCIP as QUBO solvers for Sel-Col.} 
  \def\arraystretch{1}
  \scalebox{0.81}{
    \begin{tabular}{ccrrrrrr}
    \toprule
     & & \multicolumn{3}{c}{30-sec limit} & \multicolumn{3}{c}{60-sec limit} \\
    \cmidrule(lr){3-5}\cmidrule(lr){6-8}
     & \parbox[t]{2.5cm}{\centering \phantom{o} } 
     & \parbox[t]{1.5cm}{\centering GUROBI} 
     & \parbox[t]{1.5cm}{\centering CPLEX}  
     & \parbox[t]{1.5cm}{\centering SCIP}  
     & \parbox[t]{1.5cm}{\centering GUROBI} 
     & \parbox[t]{1.5cm}{\centering CPLEX} 
     & \parbox[t]{1.5cm}{\centering SCIP}  \\
     \midrule
    \multirow{5}[0]{*}{PG} 
          & UB & 450.0 & 280.0 & 1444.8 & 239.7 & 72.3 & 1378.9 \\[0.05cm]
          & \% feas & 78.6 & 64.8 & 12.5 & 84.1 & 72.0 & 15.7 \\[0.05cm]
          & Time (sec) & 27.9 & 28.2 & 29.8 & 54.9 & 55.2 & 59.2 \\[0.05cm]
          & \% gap & 47.4 & 42.8 & 89.5 & 42.6 & 29.7 & 86.9 \\[0.05cm]
          & \% solver gap & \largeVal & \largeVal & \largeVal & \largeVal & \largeVal & \largeVal \\[0.05cm]
    \midrule
    \multirow{5}[0]{*}{ER} 
          & UB & 1710.1 & 695.0 & 2395.3 & 1507.9 & 69.0 & 2365.5 \\[0.05cm]
          & \% feas & 67.2	& 70.9	& 11.2	& 73.3	& 79.5	& 14.1 \\[0.05cm]
          & Time (sec) & 28.3 & 32.2 & 30.4 & 55.4 & 64.8 & 59.6 \\[0.05cm]
          & \% gap  & 57.9 & 57.6 & 90.3 & 56.2 & 53.3 & 87.8 \\[0.05cm]
          & \% solver gap & \largeVal & \largeVal & \largeVal & \largeVal & \largeVal & \largeVal \\[0.05cm]
    \bottomrule
    \end{tabular}%
    }
  \label{tab:selcol_summary_qubo}%
\end{table}%

%% file: conclusion.tex
\section{Conclusion}\label{conclusion}
In this work, we conduct an extensive computational study to evaluate the performance of the state-of-the-art solvers and DA on various COPs, namely, pure QUBO, quadratic assignment, quadratic cycle partition, and selective graph coloring problems, with the last two being new graph applications for DA.
For Sel-Col, we present a two-phase heuristic to reduce instance sizes and thereby improve the extent of instances that DA can handle.

Our results indicate that DA is a viable solution technology for solving problems formulated as QUBO models, usually performs comparable to or better than the state-of-the-art solvers for constrained problems, and is particularly successful in quickly yielding good-quality solutions for large instances where established solvers fail to be useful in limited times.

The number of variables that the current release of DA can handle is 8192, yet {\cfujitsu it has been recently announced that a new DA technology is developed to help handle up to one million variables \citep{da1MWeb}}.
Considering that DA demonstrates a notably robust performance despite increasing instance sizes without increased execution times, it has the potential to be a viable approach to produce high-quality solutions in a short amount of time for large-scale challenging COPs, for which even the state-of-the-art solution technologies typically struggle.

As a future research direction, the QUBO approach can be coupled with DA technology as part of different algorithms, which is a general strategy used in some existing studies, for instance \citep{liu2019modeling,shaydulin2019hybrid,shaydulin2019network,ushijimamwesigwa2020multilevel}.
Such approaches may help overcome the variable size limitation of DA and similar technologies, and also improve the overall solution performance. 
